\documentclass[a4paper]{article}
\usepackage{graphicx} 
\usepackage{amsmath, amsthm, amsfonts, amssymb, amscd, bm, enumerate}
\usepackage{verbatim}

\usepackage{color}
\usepackage{tikz}

\usepackage{url}

\renewenvironment{proof}{{\bfseries Proof:}}{\qed}
\newtheorem{theorem}{Theorem}
\newtheorem{proposition}{Proposition}
\newtheorem{corollary}{Corollary}
\newtheorem{lemma}{Lemma}
\newtheorem{definition}{Definition}
\newtheorem{assumption}{Assumption}
\newtheorem{remark}{Remark}
\DeclareMathOperator{\sign}{sign}
\DeclareMathOperator{\Tr}{Tr}

\def\B{\mathcal{B}}
\def\F{\mathcal{F}}
\def\E{\mathcal{E}}
\def\W{\mathcal{W}}

\def\L{\mathcal{L}}
\def\S{\mathcal{S}}
\def\U{\mathcal{U}}

\def\bR{\mathbb{R}}
\def\bE{\mathbb{E}}
\def\bN{\mathbb{N}}

\def\bone{\mathbf{1}}
\def\bP{\mathbb{P}}
\def\cP{\mathcal{P}}
\def\bQ{\mathbb{Q}}

\mathchardef\mhyphen="2D

\def\tT{\tilde{T}}
\def\bk{\bar{k}}
\def\<{<}
\def\>{>}

\begin{document}
\title{Ergodic BSDEs with jumps and time dependence}
\author{Samuel N. Cohen\\ University of Oxford \\ \\ Victor Fedyashov\\ University of Oxford}

\date{\today}

\maketitle

\begin{abstract}
In this paper we look at ergodic BSDEs in the case where the forward dynamics are given by the solution to a non-autonomous (time-periodic coefficients) Ornstein--Uhlenbeck-like SDE with L\'evy noise, taking values in a separable Hilbert space. We establish the existence of a unique bounded solution to an infinite horizon discounted BSDE. We then use the vanishing discount approach, together with coupling techniques, to obtain a Markovian solution to the EBSDE. We also prove uniqueness under certain growth conditions. Applications are then given, in particular to risk-averse ergodic optimal control and power plant evaluation under uncertainty.


\end{abstract}

\section{Introduction} Over the past decade, a lot of work has gone into understanding optimal control over infinite horizons. Many results for discounted problems have been obtained using techniques from classical stochastic optimal control (see, for example, Bensoussan and Lions \cite{BL_impulse}). Much less developed is the case of payoffs that value the future as much as the present, thereby being insensitive to short-term affects. One framework that has emerged is ergodic stochastic control, an area of optimal control theory that is trying to understand optimisation with an average cost criterion. Most results in this area are focused on costs which depend only on the current state of an underlying controlled Markov process, and at the linear expectation of future costs. In other words the value functional takes the form
\begin{equation}
	J(x, u) = {\lim\sup}_{T\to\infty} T^{-1} \bE^{u}\bigg[\int_0^T L(X_t, u_t) dt\bigg]
\label{opt_classical}
\end{equation}
where $X$ represents the forward dynamics, and the control $\{u_t\}_{t \geq 0}$ is an $\F_t$-predictable process taking values in a separable locally compact metric space $\U$, and $L$ is a bounded measurable cost function. It is clear that these methods are unable to deal adequately with risk-averse optimisation, since in that case a nonlinear dependence of the functional $J$ on future costs would be required. 

Since the early 90s, several papers have described the connection between Backward Stochastic Differential Equations (BSDEs), developed by Pardoux and Peng in \cite{Pardoux_Peng}, and stochastic optimal control theory (for a survey of methods see, for example, \cite{XYZ}). A strong link has also been established between BSDEs and the theory of `nonlinear expectations', as defined by Peng in \cite{Karoui_BSDE} (see Cohen \cite{Cohen_gexp} and Coquet et al. \cite{Coquet} for details). Therefore it is reasonable to expect that there exists a BSDE-based framework that would prove natural for understanding optimisation in nonlinear settings. 

One such framework is based on Ergodic BSDEs, an extension of BSDEs which takes the form 
\begin{equation}
 Y_t = Y_T+\int_t^T[f(X_{u}, Z_u)-\lambda] du - \int_t^T Z_u dW_u,
\label{EBSDE}
\end{equation}
where $\lambda \in \bR$ is a part of the solution, first introduced by Furhman, Hu and Tessitore in \cite{Hu_Banach}. Using their approach, it is relatively easy to consider nonlinear problems, for instance when the expectation $E^u$ in (\ref{opt_classical}) is replaced by a dynamically consistent nonlinear expectation (in particular, a $g$-expectation in the terminology of \cite{Karoui_BSDE}). 

The goal of present work is to extend the existing theory in two natural ways. The first generalisation is to add jumps to the diffusion setting of Furhman et al.~in \cite{Hu}. In other words, our aim is to be able to use an EBSDE--based approach to ergodic optimal control problems in the case where stochastic dynamics are given with reference to a L\'evy process. Optimal control of jump diffusions has been of great interest recently, primarily due to its possible application to network control problems and hybrid stochastic systems. From the standpoint of finance, it allows us to factor shocks into the model. The corresponding EBSDE will take the form 
\[
 Y_t = Y_T+\int_t^T[f(X_{u}, Z_u,U_{u})-\lambda] du - \int_t^T Z_u dW_u - \int_{t}^{T}\int_{H \backslash \{ 0 \} }U_{s}(x)\tilde{N}(ds,dx),  
\]
where $0 \leq t \leq T < \infty$. The second extension is to incorporate time-dependence. This will allow us to consider dynamics with seasonal components, such as business cycles. 

It is also worth noting that, since we look at EBSDEs in Markovian framework, they are related to IPDEs with nonlocal part and non-autonomous coefficients, namely
\[
\begin{cases}
-\frac{\partial}{\partial t} u(t,x) - Lu(t,x) - f(x,\nabla u(t,x)G(t), \Phi u(t,x)(\cdot)) = \lambda;  (t,x) \in \bR^+ \times H, \\
u(t + T^*,x) = u(t,x),
\end{cases}
\]
where the second-order integro-differential operator $L$ is of form
\[
	L = M + K,
\]
with
\[
	Mv(t,x) = \frac{1}{2}Tr \bigg( G(t)G^*(t) \nabla^2 v(t,x) \bigg) + \langle A(t)x + F_t(x) , \nabla v(t,x) \rangle
\]
and
\[
	Kv(t,x) = \int_{H \backslash \{0\}} \{ v(t,x + G(t)y) - v(t,x) - \langle G(t)y, \nabla u(t,x) \rangle \}\nu(dy).
\]
Derivation of this connection in finite dimensions can be found, for example, in \cite{Barles}. For equations of this type in infinite dimensional Hilbert spaces the theory is not well developed. EBSDEs provide a new way of looking at these problems. Establishing results on the connection with IPDEs is beyond the scope of present work, but it constitutes an interesting direction for future research.

The rest of the paper is organised as follows: in Section 2 we introduce the necessary notation and discuss the preliminaries; Section 3 is devoted to the results concerning solutions to the forward SDE; in Section 4 EBSDEs are introduced and main results are proven. Section 5 contains several examples of the application of EBSDEs to optimal ergodic control.

\section{Notation and general assumptions}

For the rest of the paper, let $H$ be a separable real Hilbert space with scalar product $\langle \cdot, \cdot \rangle_H$ and norm $\| \cdot \|_H$. To simplify notation we will denote them respectively $\langle \cdot, \cdot \rangle$ and $\| \cdot \|$. Since we shall be working with general separable Hilbert spaces, we will require a number of extensions of classical results. The main purpose of this section is to state them. We start with a definition of $Q$-Wiener and L\'evy processes on a general Hilbert space $H$:

\begin{definition} A stochastic process  $L = (L(t), t \geq 0)$ taking values in $H$ is called a L\'evy process if $L(0) = 0$, the process $L$ is stochastically continuous, and it has stationary, independent increments, in the sense that the law $\L(L(t)-L(s))$ depends only on the difference $t-s$. By stochastic continuity we mean that for every $\epsilon > 0$ and $t \geq 0$, $\lim_{s \to t}\bP(| L(s) - L(t) | > \epsilon) = 0$.
\end{definition}

\begin{remark}
A useful way of thinking about the L\'evy process taking values in a Hilbert space is through the series expansion, i.e. assuming that $\{ e_n \}_{n \geq 1}$ is an orthonormal basis of $H$, we have
\[
	L(t) = \sum_{n \geq 1} \langle L(t) , e_n \rangle e_n = \sum_{n \geq 1}L_n(t) e_n,
\]
where $L_n$ are real-valued c\`adlag L\'evy processes.
\label{series}
\end{remark}

\begin{definition} An $H$-valued stochastic process $\{ \W_t, t \geq 0 \}$ is called a $Q$-Wiener process if 
\begin{itemize}
\item $\W_0 = 0$,
\item $\W$ has continuous trajectories,
\item $\W$ has independent increments,
\item the law of $\W_t - \W_s$ is Gaussian with mean zero and covariance $(t-s)Q$, for all $0 \leq s \leq t$ in the sense that for any $h \in H$ and $0 \leq s \leq t$, the real-valued random variable $\langle h,  \W_t - \W_s \rangle_H$ is Gaussian with mean zero and variance $(t-s)\langle Qh,h \rangle_H$. 
\end{itemize}
\end{definition}

\noindent For a given process $\{ L_t \}_{t \geq 0}$ and a set $A \in H$ we denote by $N(t,A)$ the (random) number of `jumps of size $A$' up to time $t$, that is $N_t(A) = N(t,A) := card\{ s \in [0,t]| \Delta L_s \in A\}$. Denoting $\B(H)$ the Borel $\sigma$-algebra, we say that $A \in \B(H)$ is bounded below if $0 \notin \bar{A}$, where $\bar{A}$ denotes the closure of $A$. Proof of the following result can be found, for example, in \cite{Alru}:
\begin{proposition} If $A$ is bounded below, then $N(\cdot,A) = \{ N(t,A), t \geq 0 \}$ is a Poisson process with intensity $M(A) = \bE [N(1,A)]$. 
\end{proposition}
We remark that since we assume $H$ to be separable, it is also Polish, and therefore the space $B = H \backslash \{ 0 \}$ endowed with its Borel $\sigma$-field $\B$ is a Blackwell space. We need this since stochastic integration with respect to Poisson measures is well defined on Blackwell spaces. Following \cite{PrZab} we adopt the definition of the It\^o stochastic integral with respect to $\tilde{N}$ as an isometry, which extends the classical isometry on simple predictable processes. That is if we define
\[
	\L^2(\tilde{N}) = \bigg\{ \cP \otimes \B- \text{measurable processes } \sigma: \bE \bigg[ \int_0^t \int_B \| \sigma(s,x)\|^2 \nu(dx) ds \bigg] < \infty \bigg\}
\]
then for every $\sigma \in \L^2(\tilde{N})$ we have
\[
	\bE \bigg[ \bigg\| \int_0^t \int_B \sigma(s,x) \tilde{N}(ds,dx) \bigg\|^2 \bigg] = \bE \bigg[ \int_0^t  \int_B \| \sigma(s,x)\|^2 \nu(dx)ds \bigg].
\]
As we shall see below, any L\'evy martingale can be represented as a sum of a Wiener process and a compensated Poisson process. Therefore combining the above with the standard integration theory for Brownian motion we have a well defined stochastic integrand.

\begin{remark} It is well known that in finite dimensional spaces any L\'evy process has a c\`adl\`ag modification. However, in general this property fails in Banach spaces. But since we work with L\'evy martingales, the processes we consider can be assumed to satisfy this property (see, e.g. \cite{Zabczyk_Levy}). 
\end{remark}

\noindent The following version of the celebrated L\'evy--It\^o decomposition for an $H$-valued L\'evy process can be found, for example, in \cite{Knable}:
\begin{theorem}\textbf{(It\^o--L\'evy Decomposition)} If $L$ is an $H$-valued L\'evy process, there is a drift vector $b \in H$, a $Q$-Wiener process $\W$ on $H$ and a random measure $N$, such that $\W$ is independent of $N_{t}(A)$ for any $A$ that is bounded below, and we have
\[
	L_{t} = bt + \W(t) + \int_{||x|| < 1}x \tilde{N}(t,dx) + \int_{||x|| \geq 1}xN_{t}(dx)
\]
where $\nu$ is the L\'evy measure, and $N_{t}$ is the corresponding Poisson random measure. 
\label{Levy_Ito}
\end{theorem}
 \begin{remark} For the rest of the paper we will only be interested in the case of L\'evy martingales, and therefore the decomposition above takes the following form
 \[
 	L_{t} = \W(t) + \int_0^t \int_{B}x \tilde{N}(dt,dx).
 \]
 where $\tilde{N}(dt,dx)$ is the compensated Poisson random measure. 
 \end{remark}
\begin{assumption} Since we will mainly be dealing with square-integrable L\'evy martingales we will require the following condition to hold:
\[
	\int_{B}\|x\|^2 \nu(dx) < \infty.
\]
Given the fact that our L\'evy process is square integrable, this assumption says that there are not too many big jumps. It is not necessary in order to introduce stochastic integration with respect to L\'evy processes in a separable Hilbert space, but it will prove crucial for the coupling argument in Section  3.3.
\end{assumption}

\noindent Throughout the paper we will be repeatedly using methods involving measure changes. To that end, we need a version of the Girsanov theorem. The following is a reformulation of Theorem 15.3.10 in \cite{SCA}:

\begin{theorem} Suppose we have uniformly bounded functions $\beta: \Omega \times [0,T] \to H$ and $\gamma: B \times \Omega \times [0,T] \to \bR^+$, such that $\gamma(\cdot,\omega,t) - 1 \in L^2(\nu(dx))$ for all $(\omega,t) \in \Omega \times [0,T]$. We define
\[
	\frac{d\bQ}{d\bP} = \E \bigg(  \int_{[0,\cdot]}\beta(\omega,t)dW_t + \int_{[0,\cdot]} \int_B (\gamma(x,\omega,t) - 1) \tilde{N}(dx,dt)  \bigg)_T,
\]
where $\E$ denotes the Dol\'eans-Dade exponential. Then $\Lambda_t := \frac{d\bQ}{d\bP}\big|_{\F_t}$ is a positive square integrable martingale, and under $\bQ$
\[
	W^{\bQ} := W - \int_{[0,\cdot]} \beta(\omega,t)dt,
\]
is a Wiener process, where the integral is understood as a series (see Remark \ref{series}). The compensator of $N$ under $\bQ$ is given by 
\[
	\nu^{\bQ}(dx,dt) := \gamma(x,\omega,t)\nu(dx)dt.
\]
\label{Girsanov}
\end{theorem}

\begin{remark} In general, the assumption that $\beta$ is uniformly bounded is stronger than is necessary. However, it suffices for the purposes of this paper, since it allows us to eliminate bounded drifts by changing measure. 
\end{remark}

\section{The forward SDE}

In this section we study the properties of the `forward' process, henceforth denoted $\{X\}_{t \geq \tau}$, for some $\tau \geq 0$. Its role can be understood intuitively as a source of stochasticity in the driver of a BSDE. We first solve the dynamics of $X$ in the forward way, and then plug the obtained values into the BSDE while running it backwards. In our case, we assume that $X$ is a solution to an Ornstein--Uhlenbeck type equation driven by L\'evy noise on a separable Hilbert space $H$. We also assume that the coefficients are time periodic. This constitutes a natural way to extend the present theory and is of interest in various applications (see Chapter 6).
\subsection{Context}
We start this section by looking at a family $\{A_{t}\}_{t \geq 0}$ of linear operators on $H$ with common domain $D(A)$ dense in $H$, assuming that $A:\bR^{+} \times D(A) \to H$ generates an exponentially bounded evolution family according to the following definition (see \cite{Knable}):
\begin{definition} An exponential bounded evolution family on $H$ is a two-parameter family $\{ U(t,s) \}_{t \geq s}$ of bounded linear operators on $H$ such that 
\begin{itemize} \item $U(s,s) = I$ and $U(t,s)U(s,r) = U(t,r)$ for $r \leq s \leq t$,
\item for each $x \in H$, the map $(t,s) \to U(t,s)x$ is continuous on $s \leq t$,\text{ and}
\item there exists $M > 0$ and $\mu > 0$ such that $||U(t,s)||_{op} \leq Me^{-\mu(t-s)}$ for $s \leq t$.
\end{itemize}
\end{definition}
\begin{remark} By `generates' we mean that for $0 \leq s \leq t$ we have
\[
	\frac{d}{dt}U(t,s)x = A(t)U(t,s)x 
\]
for all $x \in H$.
\end{remark}

\begin{remark} One way of thinking about exponential bounded evolution families is as a time-dependent infinite dimensional modification of the familiar case where $A$ is a real $d \times d$ matrix, the eigenvalues of which have non-positive real parts. Then $U$ takes the form $e^{tA}$, and all conditions are satisfied.  
\end{remark}

\noindent We now consider the $H$-valued process $X$ given by the following non-autonomous mild It\^o SDE
\begin{equation}
	X(t,\tau,x) = U(t,\tau)x + \int_\tau^t U(t,s)F_s(X(s,\tau,x))ds + \int_\tau^t U(t,s)G(s)dL(s), \quad 
\label{ISDE}
\end{equation}
which is a mild version of the following Cauchy problem,
\begin{equation}
	dX_t = A(t)X_{t}dt + F_t(X_{t})dt + G(t)dL_{t}, \quad X_{\tau} = x, \quad t \geq \tau.
\label{Cauchy_problem}
\end{equation}
Conditions for existence and uniqueness of the solution will be formulated in Theorem \ref{SDE_existence}. For the rest of the paper we assume the following 
\begin{assumption} 
\begin{enumerate}[(i)]
 \item The family $A_t$ generates an exponentially bounded evolution family. Their adjoints $A^*(t)$ also have a common domain, which is dense in $H$. 
 \item $F: \bR^+ \times H \to H$ is a uniformly bounded family of measurable maps with common domain $D(F)$, which is dense in $H$. 
 \item $(\Omega,\F,\bP)$ is a complete probability space, and the pair $(W,\tilde{N})$ that comes from the It\^o--L\'evy decomposition of $L$ has the predictable representation property in the filtration $\{\F_t\}_{t\geq 0}$.
 \item $\{G_t\}_{t \geq 0}$ is a uniformly bounded family of linear operators in $L(H,H)$ with common domain $D(G)$ dense in $H$ and with bounded inverses. 
 \item The linear operator $U(t,\cdot)G(\cdot)$ is uniformly bounded in the Hilbert--Schmidt norm $\| \cdot \|_t$, defined by
 \[
 	\| S \|_t := \bigg[ \bE \bigg( \int_0^t \Tr (S_u Q S_u^*) du\bigg)  \bigg]^{\frac{1}{2}},
 \]
 where $Q$ is the covariance operator of the Wiener part of $L$.
 \item Coefficients $A(t)$, $F(t,\cdot)$ and $G(t)$ are $T^*$-- periodic for some $T^* \geq 0$, that is $A(t + T^*) = A(t)$, and similarly for $F$ and $G$.

\end{enumerate}
\label{Assumption_main}
\end{assumption}

\begin{remark} The norm $\| \cdot \|_t$ defined above allows for the following isometry:
\[
	\bE \bigg( \bigg\| \int_0^t S_td\W_t \bigg\|^2 \bigg) = \| S \|^2_t,
\]
where $\W$ is a $Q$-Wiener process and $Q$ is a trace class operator. 
\end{remark}

\noindent For the situation with autonomous coefficients, namely when $A_t = A$, $G(t) = G$ and $F(t,\cdot) = F(\cdot)$ $\forall t \geq 0$, the following theorem is a direct corollary of Theorem 9.29 in \cite{Zabczyk_Levy}:

\begin{theorem} Suppose that $A$, $F$, $G$ are time homogenous and assume that 
\begin{enumerate}[(i)]
	\item $F$ and $G$ satisfy Assumption \ref{Assumption_main},
	\item $F$ is Lipschitz-continuous.
\end{enumerate}
Then, for all $\tau \geq 0$, and any $\F_\tau$-measurable square integrable random variable $\bar{X}_\tau$ in $H$, the equation
\begin{equation}
	dX_t = (AX_t + F(X_t))dt + GdL_t, \quad X(\tau) = \bar{X}_\tau
\label{tmp1}
\end{equation}
has a unique (up to modification) mild solution with a c\`adl\`ag version. Moreover, $\forall 0 \leq \tau \leq T < \infty$, there exists $C < \infty$ such that, for all $x,y \in H$,
\begin{equation}
	\sup_{t \in [\tau,T]} \bE \|X(t,\tau,x) - X(t, \tau, y)\|^2 \leq C \| x - y \|^2.
\label{+}
\end{equation}	
\label{SDE_existence}
\end{theorem}

\begin{remark} Suppose now that $F$ is bounded and measurable, and can be approximated as a uniform limit of Lipschitz functions. Then one could adapt the proof of Theorem 10.14 in \cite{Zabczyk_Levy} to show that there exists a unique c\`adl\`ag mild solution to the equation (\ref{tmp1}). In other words, there exists an adapted $H$-valued c\`adl\`ag process $ \{X_t  \}_{t \geq \tau}$, such that the equation
\[
	X_t = e^{(t-\tau)A}x + \int_\tau^t e^{(s-\tau)A}F(X_s)ds + \int_\tau^t e^{(t-\tau)A}Gd\tilde{L}(s),
\]
is satisfied $\bP-a.s$. Moreover, the estimate (\ref{+}) still holds. 
\end{remark}
\begin{remark}
Theorem \ref{SDE_existence} can be extended to the non-autonomous case in a straightforward manner. The linear case has been treated in \cite{Knable}. For semilinear equations of the form of (\ref{ISDE}), one can prove existence by the standard fixed-point argument, and uniqueness by Gr\"onwall's lemma. Since this is not the primary interest of this work, we omit the proof. 
\end{remark}

\noindent To simplify notation, we write $U_t^s = U(t,s)$ and $U_t = U(t,0)$. Making sure that the stochastic convolutions in (\ref{ISDE}) exist in the sense of B\"ochner integral, the following result can be found, for example, in \cite{Knable}:
\begin{theorem} If $U$ is an exponentially bounded family and $G$ satisfies Assumption {\ref{Assumption_main}}, then the stochastic convolution $X_{U,G} := \int_\tau^t U(t,r)G(r)dL(r)$ exists in the following sense:
\[
	\int_\tau^t U(t,r)G(r)dL(r) = \int_\tau^t U(t,r) G(r) dW(r) + \int_\tau^t \int_{B}U(t,r)G(r)x\tilde{N}(dr,dx).
\]
\end{theorem}
\begin{definition} Whenever $f:H \to \bR$ is measurable and bounded, we call
\[
	P(s,t)[f](x) := \bE \big[ f(X(t,s,x)) \big]
\]
the two-parameter transition semigroup associated with the solution $X$ of (\ref{ISDE}). To simplify notation, in the sequel we will be particularly interested in the case $s=0$, and we write $X^x_t := X(t,0,x)$ and $\cP_t [f](x) = \bE[f(X^x_t)]$. However, all our results, including the coupling estimate, can be easily extended to the more general $P(s,t)[f](\cdot)$ case.
\end{definition}

\subsection{Coupling estimate}
The goal of this subsection is to obtain the exponential convergence of laws corresponding to two solutions of (\ref{ISDE}) with different initial conditions. We need this convergence to be uniform in the class of processes with bounded nonlinear part. In other words, our aim is to prove the following theorem:
\begin{theorem} Let $F:\bR^+ \times H \to H$ be any Lipschitz function and $\{ A_t \}_{t \geq 0}$ be fixed and generate an exponentially bounded evolution family. Then there exist constants $C > 0$ and $\rho > 0$ such that, for any bounded continuous function $\psi : H \to \bR$,
\begin{equation}
	|P(\tau,t)[\psi](x)  - P(\tau,t)[\psi](y) | \leq C (1 + ||x||^2 + ||y||^2)e^{-\rho (t-\tau)}\sup_{u \in H}||\psi(u)||
\label{Estimate_main}	
\end{equation}
where our constants $C$ and $\rho$ depend only on $\sup_{u \in H}F(u)$ and on the constant $\mu$ of the evolution family $\{ A_t \}_{t \geq 0}$.
\label{Estimate_coupl}
\end{theorem}
\noindent This estimate will be crucial in the sequel when we show the existence of a solution to an EBSDE. In our proof, we follow the derivation of Theorem 2.4 in  \cite{Hu} and Theorem 2.8 in \cite{PSXZ}. We require a number of results from the theory of coupling. A survey can be found in \cite{Lindv}. The rest of the section is organised as follows: we begin by stating the necessary facts from the theory of coupling (see \cite{Lindv} or more details). Having obtained the necessary machinery (most importantly Lemmas \ref{lemma_coupl} and \ref{coupling_bound}) we then prove Theorem \ref{Estimate_coupl}.

\begin{definition} Given two probability measures $\mu_X$ and $\mu_Y$ on measurable spaces $R_X$ and $R_Y$, a coupling is a random variable $(Z_X,Z_Y)$ taking values in the product space $R_X\times R_Y$, whose components have marginal distributions $\mu_X$ and $\mu_Y$ respectively.
\end{definition}

\begin{definition} Two processes $X$ and $Y$ are said to admit a successful coupling on $[T_1,T_2]$ when there exists $t \in [T_1,T_2]$ such that $X_t = Y_t$. 
\end{definition}

\begin{lemma} \textit{(Theorem 5.2 in \cite{Lindv})} For any two probability measures $(\mu_1,\mu_2)$ on a measurable space $(E,\E)$ there exists a coupling $(Z,Z')$ such that 
\begin{itemize} \item $\| \mu_1 - \mu_2 \|_{TV} = 2\bP(Z \neq Z')$,
			\item $Z$ and $Z'$ are independent conditional on  $\{ Z \neq Z' \}$, provided that the latter event has positive probability,
			\item $\bP(Z = Z', Z \in A) = (\mu_1 \wedge \mu_2) (A)$,
\end{itemize}
where 
\[
	\| \mu_1 - \mu_2 \|_{TV} =  \sup_{\Gamma \in \E} | \mu_1(\Gamma) - \mu_2(\Gamma) |
\]
is the standard total variation norm for measures on $(E,\E)$.
\label{coupling_max}
\end{lemma}
\begin{remark} The lemma above shows a slightly different way of thinking about couplings. Given the marginal laws, we `manually' construct random variables following them. In the process of this construction our goal is to tweak these variables in such a way as to maximise the probability of them meeting. In that case, by a coupling we mean the pair $(X,X')$ of random variables constructed. 
\end{remark}

\noindent In the sequel we will require the following auxiliary lemma, where (in principle) we couple the terminal values of solutions to (\ref{ISDE}).

\begin{lemma} Fix $\tT > 0$ and consider $X^{x,k}$ and $X^{y,k}$ the solutions (with $\tau = k\tT$) of (\ref{ISDE}) for $k \geq 0$ with initial conditions $x \in B_R(0)$ and $y \in B_R(0)$, $x \neq y$. We denote by $\mu^k_x$ and $\mu^k_y$ the respective laws of $X^{k,x}_{(k+1)\tT}$ and $X^{k,y}_{(k+1)\tT}$. Set 
\[
	Y^k_t = X_t^y + \frac{(k+1)\tT-t}{(k+1)\tT}U^{\tau}_t(x-y)
\]
and observe that the law of $Y^{k}_{(k+1)\tT}$ is $\mu^k_y$. Then, for every $p > 1$, there exists $C$ such that
\begin{equation}
	\int_H \bigg( \frac{d \mu^k_y}{d \mu^k_x}(u) \bigg)^p d\mu^k_x(u) \leq C.
\label{Cbound}
\end{equation}
\label{lemma_coupl}
\end{lemma}

\noindent \begin{proof} Given that $X^y$ satisfies (\ref{Cauchy_problem}), we immediately notice that 
\[
	dY_t = (A(t)Y_t + F_t(Y_t))dt + G(t)dL(t) - \bigg(\frac{1}{(k+1)\tT}U^{\tau}_t(x-y) + F_t(Y_t) - F_t(X_t^y)\bigg)dt.
\]
We now define 
\[
	b^*(s) = \bigg(\frac{1}{(k+1)\tT}U^{\tau}_t(x-y) + F_t(Y_t) - F_t(X_t^y)\bigg)G^{-1}(t).
\]
By assumption, $G$ is an invertible operator and  there exists $C_1 > 0$ such that $\| G^{-1} \|_{op} \leq C_1$. Given our assumptions, it is clear that 
\[
	\| b^*(s) \| \leq 2C_1 \bigg(\| F \|_{\infty} + \frac{MR}{\tT} \bigg),
\]
where $M$ comes from the definition of $U$. We define 
\[
	 \frac{d\bQ}{d\bP} = \E \bigg(  \int_{k\tT}^{(k+1)\tT} b^*(t) dW_t \bigg).
\]	
Since $b(\cdot)$ is uniformly bounded, by Theorem \ref{Girsanov}, the process $\Lambda_s := \frac{d\bQ}{d\bP} \big|_{\F_{s}}$, defined on $[k\tT, (k+1)\tT]$, is a positive square integrable martingale and $\bQ \sim \bP$. Moreover, under $\bQ$, $\tilde{L}_t = L_t - \int_0^t b^*(s)ds$ is a L\'evy process with the same triplet (in the sense of L\'evy--Khintchine representation) as $L$ under $\bP$.

\noindent It is clear that $Y_{(k+1)\tT}$ has the law $\mu^k_x$ under $\bQ$ and $\mu^k_y$ under $\bP$. We notice that
\[
	\int_H \bigg( \frac{d\mu^k_y}{d \mu^k_x}(u) \bigg)^p d\mu^k_x(u) \leq \bE \big[\Lambda_{(k+1)\tT}^p\big], 
\]
and the claim follows using that, for every $p>1$, the process $\E(\int_{k\tT}^{\cdot}p b^*(s)dW_s)$ is a true martingale. 

\end{proof}
\begin{remark} Since the coefficients in (\ref{ISDE}) depend on time, the laws of solutions with the same initial conditions on various time segments of length $\tT$ are different. However, since the bound on $b^*(\cdot)$ holds uniformly in time, the bound (\ref{Cbound}) does as well. 
\end{remark}

\noindent The following lemma can be found in \cite{Lindv}:
\begin{lemma} Let $\mu_1$ and $\mu_2$ be two equivalent probability measures on some space $E$. If there exist constants $C > 0$ and $p > 1$ such that
\[
	\int_E \bigg[ \frac{d\mu_1}{d\mu_2}(x) \bigg]^{p+1}d\mu_2(x) < C,
\]
then
\[
	\int_E \bigg( 1 \wedge \frac{d\mu_1}{d \mu_2}(x)  \bigg) d\mu_2(x) \geq \bigg[ 1 - \frac{1}{p} \bigg] \bigg( \frac{1}{pC} \bigg)^{\frac{1}{p-1}}
\]
and hence in the notation of Lemma \ref{coupling_max} we get 
\[
	\bP(Z = Z') = \big( \mu_1 \wedge \mu_2 \big) (A) \geq \bigg[ 1 - \frac{1}{p} \bigg] \bigg( \frac{1}{pC} \bigg)^{\frac{1}{p-1}}.
\]
\label{coupling_bound}
\end{lemma}

\noindent \begin{proof} (\textbf{Of Theorem \ref{Estimate_coupl}}) We concentrate on the case of $\tau = 0$ for notational simplicity. As will be clear from the proof, the result can be easily extended to the general two-parameter semigroup. We fix the initial conditions $x,y \in H$. For any two processes $X^y$ and $X^x$ with laws corresponding to the solutions of (\ref{ISDE}) with initial conditions $y$ and $x$ respectively, we denote their respective laws  $\mu_x$ and $\mu_y$. We now consider the coupled process 
\begin{equation}
	Y_t = \begin{cases} X^y_t, \quad t  < T, \\
					X^x_t, \quad t \geq T.
	\end{cases}
\label{coupled_process}
\end{equation}
where 
\[
	T = \inf\{s: X^x_s = X^y_s\}
\]
is the first meeting time of $X^x$ and $X^y$. We notice that $Y$ and $X^y$ have the same law. Now, for any bounded $\phi: H \to \bR$,
\begin{equation}
\begin{split}
	|\cP_t[\phi](x) - \cP_t[\phi](y)| &= | \bE \phi(X_t^x) - \bE \phi(Y_t)  | \\
							&= | \bE \big( [\phi(X_t^x) -  \phi(Y_t)]\bone_{\{ T > t \}} \big) | \leq 2\sup_{x \in H}|\phi(x)| \bP(T > t)
\end{split}
\label{temp}
\end{equation}
and by Markov's inequality, for any $\rho > 0$
\[
	\bP(T > t) \leq \bE\big[ e ^ {\rho T} \big] e ^ {-\rho t}.
\]
Therefore, in order to arrive at our result, we will construct $X^x$ and $X^y$, and then prove that there exist constants $\tilde{C} > 0$ and $\rho > 0$, such that 
\begin{equation}
	\bE\big[ e ^ {\rho T} \big] \leq \tilde{C} (1 + \|x\|^2 + \|y\|^2).
\label{time_bound}
\end{equation}
\begin{remark} One important thing to understand is what we mean by ``construct''. Since we are trying to prove the convergence of \textit{laws}, we do not have to work with the original solutions to our forward equation, but can instead patch together the pieces constructed on various time intervals. On each such interval $[k\tT, (k+1)\tT]$ we let $X^x_s = X(s,k\tT,X^x_{k\tT})$ with the L\'evy process $L$ in (\ref{ISDE}) replaced by $\tilde{L}$, and $X^y_s = X(s,k\tT,X^y_{k\tT})$ with $L$ replaced by $\bar{L}$, where $\tilde{L}$ and $\bar{L}$ are L\'evy processes with the same law as $L$. Since the law of the solution does not depend on the choice of the noise, $X^x$ and $X^y$ have the same laws on $[k\tT, (k+1)\tT]$ as the original solutions. 
\label{remark_construction}
\end{remark}
\noindent We proceed the following way:

\begin{itemize}
\item \textbf{(Step 1) }We start by showing that we can choose a time step $\tT > 0$ and a radius $R > 0$, such that, if we observe two independent solution processes $X^x$ and $X^y$ only at times $\{ n\tT \}_{n \in \bN}$, there is an \textit{exponential bound} on the waiting time for both $X^y_{n\tT}$ and $X^x_{n\tT}$ to enter $B_R(0)$. The independence here is understood in the the sense that we take two independent copies ($\tilde{L}$ and $\bar{L}$) of the L\'evy process $L$, as in Remark \ref{remark_construction}. 

\item \textbf{(Step 2) }Once $X^x_{k\tT}$ and $X^y_{k\tT}$ are in $B_R(0)$ for some $k \geq 0$, we lift the independence assumption and construct two solutions $X^{X^x_{k\tT}, k\tT}$ and $X^{X^x_{k\tT}, k\tT}$ to (\ref{ISDE}) on $[k\tT, (k+1)\tT]$ with initial conditions $X^x_{k\tT}$ and $X^y_{k\tT}$ respectively. We then infer that, for the constructed solutions, the probability of them meeting on $[k\tT, (k+1)\tT]$ is bounded from below uniformly in $k$.

\item \textbf{(Step 3) } We then iterate these arguments to show that the probability that the two processes we are constructing have not met decays exponentially in time.  

\end{itemize}

\noindent \textbf{Step 1: }We begin the formal derivation by showing that there exist positive constants $\mu$, $c$ and $D$ such that 
\[
	\bE \| X^x_t \|^2 \leq D(\|x\|^2e^{-2\mu t} + c).
\]


\noindent In order to proceed, we define 
\[
\begin{split}
	V_t  &= U(t,0) x + \int_0^t U(t,s)F(X_s)ds, \\
	 Z_t &= \int_0^t U(t,r) G(r) dW(r),
\end{split}
\]
and 
\[
	Q_t = \int_0^t \int_B U(t,r)G(r)x\tilde{N}(dr,dx).
\]
We notice that 
\[
\begin{split}
	\| V_t \| &\leq \| U_t x\| 	+ \bigg\| \int_0^t U_t^sF(X_s)ds \bigg\| \leq \| U_t \|_{op}\| x \| + \bar{F}\int_0^t \| U_t^s \|_{op}ds  \\
					&\leq M e^{-\mu t} \| x \| + \bar{F}\int_0^t e^{-\mu(t-s)}ds \leq M e^{-\mu t} \| x \| + \frac{\bar{F}}{\mu}.
\end{split}
\]
and thus, by using the inequality $(a + b)^2 \leq 2(a^2 + b^2)$,
\[
	\| V_t \|^2 \leq 2\big(M^2 e^{-2\mu t}\| x \|^2 + \frac{\bar{F}^2}{\mu^2} \big).
\]
Using isometries and independence of $W$ and $\tilde{N}$, we also see that
\[
\begin{split}
	\bE \| Z_t + Q_t \|^2 	 &= \bE \langle Z_t + Q_t, Z_t + Q_t \rangle_H = \bE \| Z_t \|^2 + \bE \| Q_t \|^2 \\
					& = \| U(t,\cdot) G \|^2_{t} + \int_0^t \int_B \| U_t^s G(s) x \|^2 \nu(ds) ds  \\
					& \leq C + \| G \|_{op}\int_0^t \int_B\| U_t^s \|^2_{op} \| x \|^2 \nu(dx)ds  \\
					& \leq C + \| G \|_{op}M^2\tilde{D} \leq \tilde{C},
\end{split}
\]
for some constant $\tilde{C}$, where $\tilde{D} = \int_B \| x \|^2 \nu(dx)$, $M$ comes from the definition of $U$, and $\| \cdot \|_t$ is defined as in Assumption \ref{Assumption_main}. We can now use the fact that $(a + b)^2 \leq 2(a^2 + b^2)$ again to get
\begin{equation}
	\bE \| X^x_t \|^2 \leq D(\|x\|^2e^{-2\mu t} + c)
\label{X_estimate}	
\end{equation}

\noindent for some constants $D$ and $c$. We remark that all the bounds above hold uniformly in time, that is, even though we do not have the Markov property, we still obtain that for any two solutions $X^x$ and $X^y$ of (\ref{ISDE}),
\begin{equation}
	\bE \big( \| X^x_{(k+1)\tT} \|^2 + \| X^y_{(k+1)\tT} \|^2 \big| \F_{k\tT} \big) \leq De^{-2\mu \tT}(\|X^x_{k\tT}\|^2 + \|X^y_{k\tT}\|^2) + 2Dc, \quad k \geq 0
\label{tmp1}
\end{equation}
for any fixed $\tT$. We now define, for fixed $R > 0$
\[
	A_k = \{ \| X^x_{k\tT} \|^2 + \| X^y_{k\tT} \|^2 > R \}, \quad B_k = \bigcap_{j=0}^k A_j.
\]
And by Chebyshev's inequality and (\ref{tmp1}) we obtain
\begin{equation}
	\bP (A_{k+1} \big| \F_{k\tT}) \leq \frac{De^{-2\mu\tT}}{R} ( \| X^x_{k\tT} \|^2 + \| X^y_{k\tT} \|^2) + 2\frac{Dc}{R}.
\label{tmp2}
\end{equation}
We now define the matrix
\[
	C = 
	\left(
		\begin{array}{ccc}
			D e^{-2\mu \tT} & 2Dc \\
			\frac{D}{R}e^{-2\mu\tT} & \frac{2Dc}{R}
		\end{array} 
	 \right).
\]
After multiplying (\ref{tmp1}) and (\ref{tmp2}) by $\bone_{B_k}$, taking an expectation and noticing that $\bone_{B_{k+1}} \leq \bone_{B_k}$ we have 
\[
	\left( 
		\begin{array}{ccc}
			\bE \big( \| X^x_{(k+1)\tT} \|^2 + \| X^y_{(k+1)\tT} \|^2 \big)  \bone_{B_k+1} \\
			\bP(B_{k+1})
		\end{array}
	 \right)
	 \leq C
	 \left(
	 	\begin{array}{ccc}
			\bE \big( \| X^x_{k\tT} \|^2 + \| X^y_{(k+1)\tT} \|^2 \big)\bone_{B_k} \\
			\bP(B_k)
		\end{array}
	 \right)
\]
the inequality being componentwise. Thus, iterating this procedure we arrive at
\[
	\left( 
		\begin{array}{ccc}
			\bE \big( \| X^x_{k\tT} \|^2 + \| X^y_{k\tT} \|^2 \big)  \bone_{B} \\
			\bP(B_k)
		\end{array}
	 \right)
	\leq C^k
	\left( 
		\begin{array}{ccc}
			\| x \|^2 + \| y \|^2 \\
			1
		\end{array}
	 \right),
\]
and premultiplying by the row vector $(0,1)$ on both sides we see that 
\[
	\bP(B_k) \leq (0,1) C^k
		\left( 
		\begin{array}{ccc}
			\| x \|^2 + \| y \|^2 \\
			1
		\end{array}
	 \right).
\]
The above discussion is true for any choice of $\tT$ and $R$, but now we want to obtain an exponential bound. The set of eigenvalues of $C$ is $\{ 0, \frac{2Dc}{R} + D e^{-2\mu \tT} \}$, and we need them both to be smaller than one. Therefore, we choose $R = 8Dc$ and $\tT$ such that $e^{-2\mu \tT} \leq \frac{1}{4D}$, so that $\frac{Dc}{R} + D e^{-2\mu \tT} \leq \frac{1}{2}$. Given the fact that the corresponding eigenvectors constitute a basis in $\bR^2$, the vector $(0,1)$ can be represented in the eigenvector basis, and therefore there exists a constant $\bar{C}>0$ such that 
\[
	\bP(B_k) \leq \bar{C}\left( \frac{1}{2} \right)^k(1 + \| x \|^2 + \| y \|^2).
\]
We now define the first hitting time of $B_R(0)$ on our discretised timeline as 
\[
	\tau = \inf\{k\tT: \| X^x_{k\tT} \|^2 + \| X^y_{(k+1)\tT} \|^2 \leq R, k \in \bN\},
\]
and then
\begin{equation}
	\bP(\tau \geq k\tT) \leq \bP(B_k) \leq  \bar{C}\left( \frac{1}{2} \right)^k(1 + \| x \|^2 + \| y \|^2).
\label{requrrence}
\end{equation}
Take a constant $\tilde{\beta}$ such that $\tilde{\beta}\tT < \ln 2$. Then
\[
	\bE \left( e^{\tilde{\beta}\tau} \right) = \sum_{k=0}^{\infty} e^{\tilde{\beta}k\tT}\bP(\tau = k\tT) \leq \sum_{k=0}^{\infty} e^{\tilde{\beta}k\tT}\bP(\tau \geq k\tT) \leq  \frac{\bar{C}}{1- \frac{e^{\tilde{\beta}\tT}}{2}}(1 + \| x \|^2 + \| y \|^2)
\]
and therefore there exists a constant $C_2$ such that, for every $\gamma \leq \tilde{\beta}$, 
\begin{equation}
	\bE \left( e^{\gamma \tau} \right) \leq C_2 (1 + \| x \|^2 + \| y \|^2).
\label{return_ball}
\end{equation}
The first step of the proof is now concluded. 

\noindent \textbf{Step 2: } 
We use the notation introduced in the proof of Lemma \ref{lemma_coupl}.
\begin{itemize}
\item By Lemma \ref{coupling_max}, on the interval $[k\tT,(k+1)\tT]$, there exists a pair of processes $(\tilde{X}^{x,k\tT},\tilde{Y}^{k\tT})$ with terminal time laws $\mu^k_x$ and $\mu^k_y$ respectively, such that 
\[
	\bP(\tilde{X}_{(k+1)\tT}^{x,k\tT} = \tilde{Y}^{k\tT}_{(k+1)\tT}) = \frac{1}{2}\| \mu^k_x - \mu^k_y \|_{TV}.
\]

\item We remember that we are in the case where $x,y \in B_R(0)$. Taking $p = 3$ in Lemma \ref{lemma_coupl}, and applying Lemma \ref{coupling_bound}, we know that there exists a constant $C$, such that 
\[
	\| \mu^k_x - \tilde{\mu}^k \|_{TV} \geq \frac{1}{2C},
\]
and therefore 
\[
	\bP(\tilde{X}_{(k+1)\tT}^{x,k\tT} = \tilde{Y}^k_{(k+1)\tT}) \geq \frac{1}{4C}.	
\]
\item We immediately see that the pair of processes defined by
\[
\bigg(\tilde{X}^{x,k\tT}_t,\tilde{X}^{y,k\tT}_t = \tilde{Y}^{k\tT}_t - \frac{\tT-t}{\tT}U^{k\tT}_t(x-y) \bigg)_{t \in [k\tT, (k+1)\tT]},
\]
are successfully coupled with probability bounded from below, since 
\begin{equation}
\begin{split}
	\bP(\tilde{X}^{x,k}_s &= \tilde{X}^{y,k}_s \text{ for some } s \in [k\tT, (k+1)\tT]) 
	\\ &\geq \bP(\tilde{X}^{x,k}_{(k+1)\tT} = \tilde{Y}^k_{(k+1)\tT}) \geq \frac{1}{4C}.
\end{split}
\label{bdd_below}
\end{equation}
\end{itemize}
Thus the second step is complete. 

\noindent \textbf{Step 3: }We are now ready to construct the processes $X^x$, $X^y$ we used in (\ref{coupled_process}) on individual time intervals of duration $\tT$ and then patch them all together. Assume that we have constructed $X^y$ and $X^x$ on $[0,k\tT]$. We now proceed in the following way:
\begin{itemize}
\item If $X^y_{k\tT}$ and $X^x_{k\tT}$ are in $B_R(0)$, then on $[k\tT,(k+1)\tT]$ we set $X^x_t = \tilde{X}^{X^x_{k\tT},k\tT}_t$ and $X^y_t = \tilde{X}^{X^y_{k\tT},k\tT}_t$ where $\tilde{X}$ is the maximal coupling constructed in Step 2.
\item If at least one process does not finish in the ball, then on the next timestep we set $X^x_t = \bar{X}^{X^x_{k\tT}}_t$ and $X^y_t = \bar{X}^{X^y_{k\tT}}_t$, where $\{\bar{X}^{x}_t\}_{t \geq k\tT}$ and $\{\bar{X}^{y}_t\}_{t \geq k\tT}$ are two independent solutions to (\ref{ISDE}) with $\tau = k\tT$ and initial conditions $X^x_{k\tT}$ and $X^y_{k\tT}$ respectively. 
\end{itemize}
We have thus constructed $X^x$ and $X^y$ on the entire time line. We now proceed to prove an exponential bound on their first meeting time. For that, we define a family $\{ z_k \}_{\{ k \in \bN \}}$ as follows: $z_0 = 0$ and $z_{n+1} = \inf\{k > z_n: k \in \bN, X^x_{k\tT},X^y_{k\tT} \in B_R(0)\}$. By (\ref{return_ball}) we get
\[
	\bE\big[ e ^ {\gamma z_1 \tT} \big] \leq C_2 (1 + \|x\|^2 + \|y\|^2)
\]
and thus 
\[
	\bE\big[ e ^ {\gamma (z_{n+1}-z_n) \tT} \big| \F_{z_n \tT} \big] \leq C_1 (1 + \|X^x_{z_n\tT}\|^2 + \|X^y_{z_n\tT}\|^2).	
\]
Since $e^{-\gamma z_n \tT}$ is $\F_{z_n\tT}$-measurable and $\| X^{x,y}_{z_n\tT} \| \leq R$, we get
\[
	\bE\big[ e ^ {\gamma z_{n+1}\tT }\big] \leq C_1^{n+1}C_2^n (1 + \|x\|^2 + \|y\|^2),
\]
where $C_1 = 1 + 2R^2$. We now set 
\[
	\bk = \inf\{ k :  X^x_{z_k \tT} = X^y_{z_k \tT}\}.
\]
Since $X^{x,y}_{z_k \tT} \in B_R(0)$ for every $k > 0$, we have, from (\ref{bdd_below})
\[
	\bP(\bk > k+1 | \bk > k) \leq \bigg(1 -  \frac{1}{4C} \bigg).
\]
As $\bP(\bk > k+1) = \bP(\bk > k+1 | \bk > k)\bP(\bk > k)$ we conclude that
\[
	\bP(\bk > k) \leq \bigg(1 -  \frac{1}{4C} \bigg) ^ k.
\]
We now choose $0 < \alpha < \gamma$ such that
\[
	\bigg(1 -  \frac{1}{4C} \bigg)^{1 - \alpha/\gamma} C_1^{\alpha/\gamma}C_2^{\alpha/\gamma} < 1,
\]
and then, using H\"older's inequality, we see that
\[
\begin{split}
	\bE(e^{\alpha z_{\bk} \tT}) &= \bE \big( \bE (e^{\alpha z_{\bk} \tT} | \bk) \big)  
	 \leq \sum_{k \geq 0} \bE \big[ e^{\alpha z_k \tT} \bone_{\bk = k}\big] \\
	&\leq \sum_{k \geq 0} (\bP(\bk = k))^{1 - \alpha/\gamma} (\bE e^{\alpha z_k \tT})^{\alpha/\gamma}  \\
	& \leq \sum_{k \geq 0} (\bP(\bk > k-1))^{1 - \alpha/\gamma} (\bE e^{\alpha z_k \tT})^{\alpha/\gamma}  \\
	& \leq \sum_{k \geq 0} \bigg(1 -  \frac{1}{4C} \bigg) ^ {(k-1)(1 - \alpha/\gamma)} (C_1^{k}C_2^{k-1} (1 + \|x\|^2 + \|y\|^2))^{\alpha/\gamma} \\
	& \leq C_3  (1 + \|x\|^2 + \|y\|^2)
\end{split}
\]
for some constant $C_3$. For each $\rho \leq \alpha$ we get
\[
	\bE\big[ e ^ {\rho T} \big] \leq \bE\big[ e ^ {\rho (z_{\bk} + 1)\tT} \big] \leq  \tilde{C}  (1 + \|x\|^2 + \|y\|^2)
\]
where $\tilde{C} = C_3 e^{\rho \tT}$.
\end{proof}

\begin{lemma} The estimate (\ref{Estimate_main}) can be extended to the case where $F$ is bounded and measurable, and there exists a uniformly bounded sequence of Lipschitz (in the second argument) functions $\{ F_n \}_{n \geq 1}$ such that 
\[
	\lim_{n} F_n(t,x) = F(t,x), \quad \forall x \in H, t \geq 0.
\]
\label{Lemma_tmp}
\end{lemma}
\noindent \begin{proof} The proof uses standard Girsanov arguments and is identical to Corollary 2.5 in \cite{Hu}.

\end{proof}

\begin{remark} The reason Lemma \ref{Lemma_tmp} is necessary is because in order to construct a solution to the EBSDE in the sequel, we will first have to change measure. From Girsanov's theorem, we know that under the new measure the forward process will have additional bounded drift. We therefore need to ensure that the estimate (\ref{Estimate_main}) still holds. 
\end{remark}

\subsection{Recurrence}
This section is devoted to proving that under certain assumptions the forward process (\ref{ISDE}) eventually enters any open ball in $H$ with probability one. We establish this for the case of time periodic coefficients and L\'evy noise with non-trivial diffusion component. This is a natural extension of existing theory and interesting in its own right. In the sequel we will need a slightly weaker property, namely the eventual return to any open ball around zero, in order to prove the uniqueness of the Markovian solution to an EBSDE. We start by formulating an additional assumption:
\begin{assumption} For notational simplicity suppose that in (\ref{ISDE}) $\tau = 0$. Then we assume that the process $Z_A(t)$ defined by
\[
	Z_A(t) = \int_0^t U_s G(s)dL_s.
\]
spans the entire space $H$, that is, $\bP(Z_A(t) \in V) > 0$ for all $t > 0$ and any open $V \in H$.
\label{span_H}
\end{assumption}
\begin{remark} This assumption may seem overly restrictive, as one can think of many L\'evy processes that do not span the entire space. For example, the case when one-dimensional components $\{L_n(t)\}$ are supported on the integers. Even in a more general case, one could think of a L\'evy process $L(t)$ supported on a subspace. However, since we focus our attention on the case where $G(s)$ is invertible for every $s \geq 0$, and $L$ has a non-trivial diffusion component, the assumption is reasonable.  
\end{remark}

\begin{lemma} If the process $Z_A(t)$ satisfies Assumption \ref{span_H} for all $t > 0$, then process $X^x_t$ satisfying (\ref{ISDE}) is irreducible. In other words
\[
	\bP(X^x_t \in B_{\epsilon}(z)) > 0
\]
for any $t > 0$, $z \in H$, $\epsilon > 0$.
\label{irreducibility}
\end{lemma}
\begin{remark} Here, and in the sequel, we denote by $B_R(x)$ the open ball of radius $R$ around some $x \in H$.
\end{remark}
\noindent \begin{proof} We follow the proof of Proposition 3.3 in \cite{PSXZ}. We fix $T > 0$, $y \in H$, $\epsilon > 0$. For the rest of the proof we also denote $X_t = X^x_t$. Then
\[
	X_{t+a} = U_{t+a}^t X_t + \int_t^{t+a} U_{t+a}^s F_s(X_s)ds + \int_t^{t+a}U_{t+a}^s G(s) dL_s.
\]
Let $z$ be any element in the support of the distribution of the random variable $U_{t+a}^s X_t$. Then, by definition, the event 
\[
	B = \{|U_{t+a}^s X_t - z| < \epsilon/3\}
\]
is of positive probability. Since $\|  F \|_\infty = \sup_{t \geq 0, x \in H} F_t(x) < \infty$, using the definition of $U$, we have
\[
\begin{split}
	\bigg | \int_t^{t+a} U_{t+a}^s F_s(X_s)ds \bigg | &\leq  \int_t^{t+a} \| U^s_{t+a} \|_{op} \| F \|_{\infty}ds  \\
	 &\leq  M \|  F \|_\infty \int_t^{t+a}e^{-\mu (t+a-s)}ds \\
	 &\leq ca,
\end{split}
\]
for some $c > 0$. We then write
\begin{equation}
	X_{t+a} - y = (U_{t+a}^t X_t - z) + \int_t^{t+a} U_{t+a}^s F_s(X_s)ds + \bigg(\int_t^{t+a}U_{t+a}^s G(s) dL_s - y + z\bigg).
\label{*}
\end{equation}
The event 
\[
	C = \bigg\{ \bigg| \int_t^{t+a}U_{t+a}^s G(s) dL_s - y + z \bigg| < \epsilon / 3 \bigg\}
\]
is of positive probability by Assumption \ref{span_H}. Since $X_t$ and the increments of $L$ on $[t,t+a]$ are independent, so are the events $B$ and $C$. Therefore $B \cap C$ has positive probability. Given (\ref{*}), we have shown that
\[
	|X_{t+a} - y| \leq \epsilon/3 + ca + \epsilon/3
\]
with positive probability on $B \cap C$. We now choose $a$ so that $ca < \epsilon/3$ and $T-a \geq 0$. Setting $t = T-a$, we obtain 
\[
	\bP(|X_T - y| \leq \epsilon) \geq \bP(B \cap C) > 0,
\]
which is the result.

\end{proof}

\noindent In order to proceed, we need a few results concerning the invariant measure for the solution to the equation (\ref{ISDE}). We begin by considering the linear problem
\[
	dX_t = A(t)dt + G(t)dL_t, \quad X_\tau = x,
\]
which can be reduced to the autonomous case by the standard technique of enlarging the state space, i.e. by considering the evolution of the vector $(X,y) \in H \times \bR_{+}$ given by
\[
\begin{cases}
	dX_t = A(y(t))X(t) + G(y(t))dL_t \quad  &X(0) = x \\
	dy(t) = dt &y(0) = \tau
\end{cases}
\]
Following \cite{Knable} we define a one-parameter semigroup as
\[
	P_s u(t,x) := P(t,t + s)u(t+s,\cdot)(x)
\]
meaning that we apply the two-parameter semigroup to $u$ as a function of $x$ only. It is clear from the definition that $P_\tau$ is a Markovian semigroup, which gives us the opportunity to use the powerful existing theory. In order to establish existence and uniqueness of the invariant measure, we need to define the corresponding ``periodic" $L^2$-space on which the semigroup is a contraction. We denote 
\[
\begin{split}
	L_*^2(\nu) := \bigg\{ &f:\bR \times H \to \bR \text{ measurable } : f(t+T^*,x) = f(t,x) \quad \nu -a.e. \\
	& \text{ and } \int_{[0,T^*] \times H} |f(y)|^2\nu(dy) < \infty \bigg\}.
\end{split}
\]
for some measure $\nu$. It is clear that $L_*^2$ is a Hilbert space. The following result was established in \cite{Knable}.
\begin{proposition} There exists a unique invariant measure for the semigroup $P$. In other words for every bounded measurable function $u$ such that $u(t+T^*,x) = u(t,x)$ for each $t >0$ and $x \in H$ we have:
\[
	\int_{[0,T^*]\times H} P_s u(t,x) \nu(dt,dx) = \int_{[0,T^*] \times H}u(t,x)\nu(dt,dx).
\]
Furthermore, on $L_*^2(\nu)$ the semigroup $P_s$ is a contraction.	
\label{invariant_existence}
\end{proposition}

\noindent We deduce that there also exists a unique invariant measure $\mu$ corresponding to the original semilinear problem
\[
	dX_t = A(t)dt + F_t(X_t)dt + G(t)dL_t, \quad X_\tau = x.
\]
This can easily be shown by a change of measure to reduce to the linear case. We leave details to the reader. We are now ready to prove the main result of this section, namely the recurrence of the forward process $\{ X_t \}_{t \geq \tau}$. We present two proofs: one is applicable for the case of $\dim H < \infty$, and is elementary, in the sense that it does not rely on the existence of the invariant measure. The second one deals with the case of $\dim H = \infty$. 

\begin{theorem} For any $x_0,x \in H$, $s \geq 0$ and for any fixed $\epsilon > 0$,  we define $\tau := \inf\{t \geq s: X^x_t \in B_\epsilon(x_0)\}$. Then $\bP(\tau > T) \to 0$ as $T \to \infty$. 
\label{recurrence}
\end{theorem}

\noindent \begin{proof} \textbf{(Intuition, $\dim H < \infty$)} From Step 1 of the proof of Theorem \ref{Estimate_coupl} we know that we can find a radius $R > 0$, such that the probability that our process returns to the ball $\bar{B}_R(0)$ is not trivial. We then discretise time with a step $\tT$. We know that the discretised process will return to $\bar{B}_R(0)$ infinitely often, and by Lemma \ref{irreducibility} the probability of the jump from $\bar{B}_R(0)$ to any open ball is bounded from below. We then invoke a Borel--Cantelli type argument to demonstrate the claim. 

\noindent \textbf{(Formal proof, $\dim H < \infty$)} We start by introducing a family of events $\{ E_n \}_{n \geq 1}$ as
\[
	E_n = \{ \text{there exists } k = 1 \dots n : X^x_{k\tT} \in B_\epsilon(x_0) \}, 
\]
and immediately notice that
\[
	\bP (E_n \big| \bar{E}_{n-1}) =  \bP( \{ X(n\tT, (n-1)\tT, X^x_{(n-1)\tT}) \in B_\epsilon(x_0) \}),
\]
where $\bar{E}$ denotes the complement of $E$ and $X(t,s,x)$ is the value at time $t$ of the solution to (\ref{ISDE}) starting at time $\tau = s$ with $X_\tau = x$.  Therefore,
\[
\begin{split}
	\bP(E_n | \bar{E}_{n-1}) &= \bP \big( X(n\tT, (n-1)\tT, X^x_{(n-1)\tT}) \in B_\epsilon(x_0) \big) \\
					&\geq \bP \big( X^x_{(n-1)\tT} \in B_R(0), X(n\tT, (n-1)\tT, X^x_{(n-1)\tT}) \in B_\epsilon(x_0) \big) \\
					     &= \bP	\big(X^x_{(n-1)\tT} \in B_R(0)) \bP( X(n\tT, (n-1)\tT, X^x_{(n-1)\tT}) \in B_\epsilon(x_0) \big).
\end{split}
\]
Since coefficients in (\ref{ISDE}) are $T^*$-periodic and $\bar{B}_R(0)$ is compact (and therefore $[0,T^*] \times \bar{B}_R(0)$ is compact), and given the stability of solutions to (\ref{ISDE}) with respect to the initial value (as stated in (\ref{+})), there exists $\delta >0$ such that
\[
	\bP \bigg(X(n\tT,(n-1)\tT,X^x_{n\tT}) \in B_\epsilon(x_0) \bigg| X^x_{n\tT} \in \bar{B}_R(0) \bigg) > \delta.
\] 
Hence 
\[
	\sum_{n \geq 1} \bP(E_n | \bar{E}_{n-1}) \geq \delta \sum_{n \geq 1} \bP(X^x_{n\tT} \in \bar{B}_R(0)).
\]
In ``Step 1" of the the proof of Theorem \ref{Estimate_coupl} we showed that 
\[
	\bE \| X^x_t \|^2 \leq L(\|x\|^2e^{-2\mu t} + c),
\]
and therefore by Markov's inequality we have
\[
	\bP(\| X^x_t \|^2 > R) \leq \frac{L}{R}(\|x\|^2e^{-2\mu t} + c).
\]
It is clear that we can choose $R$ so that 
\[
	1 - \bP(\| X^x_{k\tT} \|^2 > R) \geq 1/k
\]
for all $k \geq 1$. Therefore 
\[
	\sum_{n \geq 1} \bP(E_n | \bar{E}_{n-1}) \geq \delta \sum_n 1/n = \infty,
\]
and thus by the counterpart of the Borel--Cantelli Lemma (see \cite{Bruss}), we conclude that
\[
	\bP(\tau < \infty) = \bP(\cup_{n}E_n) = 1,
\]
concluding the proof.
\end{proof}

\noindent \begin{proof} (\textbf{$\dim H = \infty$}) Let the function $\psi: \to \bR$ be bounded and continuous. From Theorem \ref{Estimate_main} we know that for any $x,y \in H$ and $0 \leq t \leq t'$ we have
\begin{equation}
\begin{split}
	|P_s[\psi](t,x)  - P_s[\psi](t',y) | &= |P_s[\psi](t',X^{t,x}_{t'})  - P_s[\psi](t',y) | \\
							&\leq C (1 + ||X^{t,x}_{t'}||^2 + ||y||^2)e^{-\rho (s-t')}\sup_{u \in H}||\psi(u)||.
\end{split}
\label{e1}
\end{equation}
Using (\ref{e1}) and the fact that there exists a unique invariant measure $\nu$ for the semigroup $P_t$, one can show (see, e.g. \cite{G.DaPrato1996}) that $\nu$ is exponentially mixing. In other words,
\[
	P_s[\psi](t,x) \to \nu(\psi) = \int_{[0,T^*] \times H} \psi(t,x) \nu(dt,dx).
\]
We now set $\psi(t,x) = \bone_{x \in A}$ for some open set $A \subset H$. Then $P_s[\psi](t,x) = \bP (X^{t,x}_s \in A)$. By Theorem \ref{irreducibility} we know that for all $0 \leq t \leq s$, $x \in H$ we have $\bP (X^{t,x}_s \in A) > 0$. Therefore 
\[
	\nu([0,T^*] \times A) = \int_{[0,T^*]\times H} \bP (X^{t,x}_s \in A) \nu(dt,dx) = \delta_A > 0
\]
for some constant $\delta_A$. Setting $t = 0$ and $A = B_\epsilon(x_0)$ we have 
\[
	\liminf_{t\to\infty}\bP(X^x_t \in B_\epsilon(x_0)) = \nu([0,T^*]\times B_\epsilon(x_0)) = \delta_{\epsilon} > 0
\]
and thus, by Proposition 3.4.5 in \cite{PRATO}, the claim follows. 

\end{proof}

\section{Backwards SDEs}
We now move from the `forward' process $X$ to consider the `backwards' part of our problem. This section is organised as follows: we start by introducing the class of discounted BSDEs in infinite horizon and proving that they admit a bounded solution. Then we use the coupling estimate obtained in the previous section to prove existence of a solution to our EBSDE. The next subsection is devoted to the uniqueness of the Markovian solution. We conclude by providing an alternative representation for the solution. Similarly to \cite{Royer_jumps} we impose certain assumptions on the driver of our BSDE.

\begin{definition} Henceforth we assume that the driver of a BSDE with jumps is a measurable function $f:\Omega \times \bR^+ \times \bR \times H \times \L^2(B,\B,\nu) \to \bR$.
\end{definition}

\begin{assumption} For all $T$ we have the following conditions on our driver $f(\omega,t,y,z,u)$:
\begin{itemize}
\item $f$ is predictable in $(\omega,t)$.
\item $f$ is continuous w.r.t y and there exists an $\bR^+$-valued process $(\phi_{t})_{0\leq t \leq T}$ such that $\bE\left(  \int_{0}^{T} \phi^{2}_{s}ds \right) < \infty$ and 
\[	
	|f(\omega,t,y,z,u)| \leq \phi_{t} + K \left(  |y| + ||z|| + \int_{B} |u(v)|^{2} \nu(dv) \right)^{1/2}
\]
\item $f$ is ``monotonic" w.r.t y, that is $\exists \alpha \in \bR$ such that $ \forall t \geq 0, \forall y,y' \in \bR, \forall z \in H, \forall u \in \L^2(B,\B,\nu)$
\[
	(y-y')(f(\omega,t,y,z,u) - f(\omega,t,y',z,u)) \leq \alpha |y - y'|^{2} \quad \bP-a.s.
\]
\item $f$ is Lipschitz w.r.t. z and u. In particular $\exists K \geq 0 : \forall t \in [0,T], \forall y \in \bR, \forall z,z' \in H, \forall u,u' \in \L^2(B,\B,\nu)$
\[
	|f(\omega,t,y,z,u)-f(\omega,t,y,z',u')| \leq K|| z-z' || + K \left(  \int_{B}|u(v)-u'(v)|^2 \nu(dv) \right)^{1/2}
\]

\end{itemize}
\label{ass:driver}
\end{assumption}

\noindent In order to have a comparison theorem, we make the following further assumption.

\begin{assumption}\label{a:generator}There exists $-1<C_{1}\leq 0$ and $C_{2}\geq0$ such that
\[
	\forall x \in H, \quad \forall z \in H^*, \quad \forall u,u' \in \L^{2}(B,\B,\nu,\bR)
\]
we have
\[
	f(\omega,t,y,z,u) - f(\omega,t,y,z,u') \leq \int_{B}(u(v)-u'(v))\gamma_t^{\omega,z,u,u'}(v)\nu(dv),
\]
where $\gamma^{\omega,t,z,u,u'}:\Omega \times B \to \bR$ is measurable (predictable $\times$ Borel) in all arguments and satisfies 
\[
	C_{1}(1 \wedge ||x||) \leq \gamma^{\omega,t,z,u,u'}(x) \leq C_{2}(1 \wedge ||x||)
\]
for all $v \in B$.
\label{ass:gamma}
\end{assumption}

\noindent The following existence theorem for finite horizon BSDEs with jumps can be found in \cite{Royer_jumps}. In that paper the case of finite-dimensional Brownian motion is considered. The extension to the infinite dimensional case where $W$ is a $Q$-Wiener processes is immediate (for details see \cite{SCA}). 
\begin{theorem} Under Assumption \ref{ass:driver}, there exists a unique solution $(Y,Z,U) \in (\S^2 \times \L^2(W) \times \L^2(\nu))$, for any terminal condition $\eta \in \L^2(\F_T)$, to the equation
\[
 Y_t = \eta + \int_t^T f(\omega,u,Y_u, Z_u,U_{u}) du - \int_t^T Z_u^* dW_u - \int_{t}^{T}\int_{B}U_{s}(x)\tilde{N}(ds,dx).
\]
\end{theorem}

\begin{lemma} \label{lm:gammaProcess} For every $Y,Z,U,U'$ under Assumption \ref{a:generator} there exists a process $\gamma_{t} = \gamma^{\omega,Y_t,Z_{t},U_{t},U'_{t}}$ such that 
\begin{equation}
	f(\omega,t,Y_t,Z_{t},U_{t}) - f(\omega,t,Y_t,Z_{t},U'_{t})= \int_{B}(U(v)-U'(v))\gamma_{t}(v)\nu(dv)
\label{change_jump}
\end{equation}
\end{lemma}
\noindent \textbf{Proof:} We first notice that $\forall t > 0, \forall \omega \in \Omega,  \forall z \in H^*, \forall u,u' \in \L^{2}(B,\B,\nu,\bR)$ there exist $\gamma_{1,t}^{\omega,z,u,u'}(v)$ and $\gamma_{2,t}^{\omega,z,u,u'}(v)$, satisfying $C_{1}(1 \wedge ||v||) \leq \gamma_{i,t}(v) \leq C_{2}(1 \wedge ||v||)$ for $i=1\dots 2$, such that 
\[
	f(\omega,t,y,z,u) - f(\omega,t,y,z,u') \leq \int_{B}(u(v)-u'(v))\gamma_{1_t}^{\omega,y,z,u,u'}(v)\nu(dv)
\]
and
\[
	f(\omega,t,y,z,u) - f(\omega,t,y,z,u') \geq \int_{B}(u(v)-u'(v))\gamma_{2,t}^{\omega,y,z,u,u'}(v)\nu(dv).
\]
Then there exists $\alpha_t = \alpha(t,\omega,y,z,u,u')$ such that 
\[
\begin{split}
	f(\omega,t,y,z,u) &- f(\omega,t,y,z,u') \\ 
		&= \int_{B}(u(v)-u'(v))(\alpha_t  \gamma_{1,t}^{\omega,z,y,u,u'}(v) + (1-\alpha_t)\gamma_{2,t}^{x,z,y,u,u'}(v))\nu(dv)
\end{split}
\]
and we immediately see that 
\[
	\alpha_t = \frac{f(\omega,t,y,z,u) - f(\omega,t,y,z,u') - \int_{B}(u(v)-u'(v))\gamma_{2,t}^{\omega,z,y,u,u'}(v)\nu(dv)}{\int_{B}(u(v)-u'(v))(\gamma_{1,t}^{\omega,y,z,u,u'}(v)-\gamma_{2,t}^{\omega,y,z,u,u'}(v))\nu(dv)} \in [0,1],
\]
noticing that if the denominator is zero then $\alpha_t = 1$ satisfies the claim. Now for each $s \in [0,t]$ and $v \in B$ we can explicitly define 
\[
	\gamma_{t}(v)= \alpha(t,\omega) \gamma_{1_t}^{\omega,Y_t,Z_{t},U_{t},U'_{t}}(v) + (1 -  \alpha(t,\omega)) \gamma_{2,t}^{\omega,Y_t,Z_{t},U_{t},U'_{t}}(v)
\]
where $ \alpha(t,\omega) = \alpha (t,\omega,Y_t,Z_{t},U_{t},U'_{t})$, and it is clear that $\gamma_{t}$ satisfies (\ref{change_jump}). 

\qed

\subsection{Infinite horizon BSDEs}

In this section we show that there exists a unique bounded solution to the infinite-horizon BSDE with discounting, that is the equation 

\begin{equation}\label{eq:discBSDE}
Y_T=Y_t - \int_{t}^{T} (-\alpha Y_{u} + f(\omega, u, Z_u, U_u)) du + \int_{t}^{T} Z_u dW_u  + \int_{t}^{T}\int_{B}U_{s}(x)\tilde{N}(ds,dx),  
\end{equation}
which will prove crucial to the study of Ergodic BSDEs in the next section. In order to proceed we will require Tanaka's formula for general semimartingales. The following version can be found, for example, in \cite{SCA}. Here we use the convention that $\sign(x)=x/|x|$ for $x\neq 0$ and $\sign(0)=0$.
\begin{lemma}\textbf{(Tanaka's formula)} Let $X$ be a semimartingale and $a \in \bR$. Then there exists a continuous increasing local time process $L^{a}$ and a pure jump process $L^{X,a}$, with $L^{a}(0) = 0$ (unique $\bP-a.s.$), such that $X$ allows the following representation: 
\[
	d|X_{t}-a| = \sign(X_{t-}-a)dX_{t} + dL^{a}_{t} + \Delta L^{X,a}_{t},
\]
where 
\[
	\Delta L^{X,a}_{t} = |X_{t}-a| - |X_{t-}-a| - \sign(X_{t-}-a)\Delta X_{t}
\]
is a `local-time' jump process.
\end{lemma}
\begin{remark} If we consider the above process $\Delta L^X$, we notice that
\[
\begin{split}
	\Delta L^{X,a}_t   &= |X_{t}-a| - |X_{t-}-a| - \sign(X_{t-}-a)\Delta (X_{t}-a) \\
			      	&= |X_{t}-a| - \sign(X_{t-}-a)(X_{t}-a) \\ 
				&= |X_{t}-a|(1 - \sign((X_{t-}-a)(X_{t}-a))) \\
				& \geq 0.	
\end{split}
\]
\label{rm:Tanaka}
\end{remark}

\begin{theorem}\label{thm:discountedBSDEexist}
Let $\alpha>0$ and $f:\Omega\times\bR^+\times H^* \times \bR \to \bR$ be such that 
\begin{itemize} 
\item $f$ satisfies Assumptions \ref{ass:driver} and \ref{ass:gamma}
\item $|f(w,t,0,0)|$ is uniformly bounded by $C \in \bR$
\end{itemize}
Then there exists an adapted solution $(Y, Z,U)$, with $Y$ c\`adl\`ag and $Z\in \L^2(W)$, $U \in \L^2(\tilde{N})$to the infinite horizon equation (\ref{eq:discBSDE}) for all $0\leq t\leq T<\infty$, satisfying $|Y_t|\leq C/\alpha$, and this solution is unique among bounded adapted solutions.

Furthermore, if $(Y^T, Z^T,U^T)$ denotes the (unique) adapted square integrable solution to
\begin{equation}\label{eq:ThorizonBSDE}
 Y^T_t = \int_{t}^{T} (-\alpha Y_{u}^T + f(\omega, u, Z_u^T, U^T_u)) du - \int_{t}^{T} (Z_u^T)dW_u - \int_{t}^{T}\int_{B}U_{s}(x)\tilde{N}(ds,dx)
\end{equation}
then $\lim_{T\to\infty} Y^T_t = Y_t$ a.s., uniformly on compact sets in $t$.
\end{theorem}
\noindent \begin{proof} We start by proving that if a bounded solution exists, it is unique. Suppose we have two bounded solutions $(Y,Z,U)$ and $(Y',Z',U')$ to (\ref{eq:discBSDE}). We denote $\delta Y := Y - Y'$, $\delta Z := Z - Z'$ and $\delta U = U - U'$. We also denote 
\[
	\alpha_s := \begin{cases}
			   	\frac{f(\omega,s,Z_s,U_s)-f(\omega,s,Z'_s,U_s)}{||\delta Z||^2}\delta Z & \text{if } Z \neq Z', \\
				0 & \text{otherwise.}	
			  \end{cases}
\]
Now define $M_t = \int_{0}^{t}\alpha_{s}dW_s + \int_{0}^{t}\int_{B}\gamma_s(x)\tilde{N}(ds,dx)$, where $\gamma = \gamma^{\omega,t,Z',U,U'}$ is defined as in Lemma \ref{lm:gammaProcess}. Then we can use Theorem \ref{Girsanov} to show that there exists a probability measure $\bQ \sim \bP$ such that under $\bQ$ the process
\[
\begin{split}
	K_t = \int_{t}^{T} (f(\omega, u, Z_u, U_u)&-f(\omega, u, Z'_u, U'_u)) du \\
	&+ \int_{t}^{T} \delta Z_u^*dW_u  + \int_{t}^{T}\int_{B}\delta U_{s}(x)\tilde{N}(ds,dx)
\end{split}
\]
is a martingale. We now apply Tanaka's formula and Remark \ref{rm:Tanaka} to see that for all $s \leq t \leq T$ we have 
\[
	E_{\bQ}[e^{-\alpha t}|\delta Y_t| - e^{-\alpha s}|\delta Y_s|\,|\F_s] \geq 0,
\]
and hence 
\[
	|\delta Y_s| \leq e^{-\alpha (t-s)}E_{\bQ}[|\delta Y_t|\,|\F_s] \leq e^{-\alpha (t-s)} C,
\] 
for $C$ a bound on $|\delta Y_t|$. This bound is independent of $T$ and collapses as $t\to\infty$. Hence $|\delta Y_s|=0$, from which we see $Y_s=Y'_s$ a.s. for every $s$, and hence $Y=Y'$ up to indistinguishability as $Y$ and $Y'$ are c\`adl\`ag.

We now show that a bounded solution exists. We first notice that there indeed exists a unique solution to the $T$-horizon BSDE (\ref{eq:ThorizonBSDE}). In order to see this, by a standard comparison argument it suffices to check that our new driver, namely $F(\omega,t,y,z,u) := -\alpha y + f(\omega, t, z, u)$ satisfies Assumption \ref{ass:driver}, provided that $f$ does. This is clear given that our additional term does not depend on $(z,u)$, is continuous and monotonic. We denote the solution as $(Y^T,Z^T,U^T)$. We now prove that $Y^T$ is bounded. Similarly to above, we define
\[
	\beta_s := \begin{cases}
			   	\frac{f(\omega,s,Z^T_s,U^T_s)-f(\omega,s,0,U^T_s)}{||Z^T||^2} Z^T_s & \text{if } Z^T_s \neq 0, \\
				0 & \text{otherwise.}	
			  \end{cases}
\]
and denote $\tilde{M}_t = \int_{0}^{t}\beta_{s}dW_s + \int_{0}^{t}\int_{B}\tilde{\gamma}_s(x)\tilde{N}(ds,dx)$, where $\tilde{\gamma} = \gamma^{\omega,0,U,0}$ is defined as in Lemma \ref{lm:gammaProcess}. Then, as above, there exists a probability measure $\tilde{\bQ} \sim \bP$, under which the process
\[
	\tilde{K}_t = \int_{t}^{T} (f(\omega, u, Z_u, U_u)-f(\omega, u, 0, 0)) du + \int_{t}^{T}  Z_u^*dW_u  + \int_{t}^{T}\int_{B}U_{s}(x)\tilde{N}(ds,dx)
\]
is a martingale, and therefore applying Tanaka's formula and It\^o's formula to $e^{-\alpha t}|Y^T_t|$ we see that
\begin{equation}\label{eq:boundsoln}
 |Y^T_t| \leq e^{\alpha t}E_{\tilde{\bQ}}\Big[\int_t^T e^{-\alpha u}|f(\omega, u, 0,0)|du\Big|\F_t\Big]\leq C/\alpha
\end{equation}
where $C$ is the bound on $|f(\omega,t,0,0)|$. Thus $Y^T$ is uniformly bounded. We now show that $Y^T$ forms a Cauchy sequence in $T$ uniformly on compacts in t. For every $T' \geq T$ we define 
\[
	\rho_s := \begin{cases}
			   	\frac{f(\omega,s,Z^T_s,U^T_s)-f(\omega,s,Z^{T'},U^{T}_s)}{||Z^T-Z^{T'}||^2} (Z^T_s - Z^{T'}_s) & \text{if } Z^T_s \neq Z^{T'}_s, \\
				0 & \text{otherwise,}	
			  \end{cases}
\]
and denote $\bar{M}_t = \int_{0}^{t}\beta_{s}dW_s + \int_{0}^{t}\int_{B}\bar{\gamma}_s(x)\tilde{N}(ds,dx)$, where $\bar{\gamma} = \gamma^{\omega,Z^{T'},U^{T},U^{T'}}$ is defined as in Lemma \ref{lm:gammaProcess}. As above, applying Tanaka's formula and inequality (\ref{eq:boundsoln}), we observe 
\begin{equation}
	|Y^T_t- Y^{T'}_t| \leq e^{-\alpha(T-t)} E_{\bar{\bQ}}[|Y^{T'}_t - Y^T_t|\,|\F_t] \leq 2Ce^{-\alpha(T-t)}/\alpha, 
\label{conv}
\end{equation}
where $\bar{\bQ}$ is defined in a similar way to $\bQ$. Hence we see that $Y^T_t$ is a Cauchy sequence in $T$, therefore the limit exists, and we denote it $Y_t$. The bound established in inequality (\ref{eq:boundsoln}) also holds for $Y_t$, and convergence uniformly on compacts is clear from (\ref{eq:boundsoln}). We now show that $Z^T_t$ and $U^T_t$ are Cauchy sequences. We denote $\tilde{Z}_t = Z^T_t - Z_t^{T'}$ and $\tilde{U}_t = U^T_t - U_t^{T'}$. Apply It\^o's formula to $\big(\tilde{Y_t}\big)^2$, where $\tilde{Y_t} := Y^T_t-Y^{T'}_t$. Then, after standard calculations under $\bar{\bQ}$, we see that, for each $t < T$,
\[
	\tilde{Y}_t^2 = \tilde{Y}_0^2 + \bE_{\bar{\bQ}} \left( \int_0^t \int_{B} |\tilde{U}_s(v)|^2 \nu(dv)ds\right) + \bE_{\bar{\bQ}} \bigg( \int_0^t ||\tilde{Z}_s||^2 ds \bigg) + 2 \alpha \bE_{\bar{\bQ}} \bigg( \int_0^t \tilde{Y}_s^2 ds \bigg).
\]
Given (\ref{conv}) our claim follows. Therefore the limit as $T \to \infty$ exists for sequences $\{Z^T_t\}$ and $\{U^T_t\}$. Taking $Z$ and $U$ as their respective limits, we have our desired solution $(Y,Z,U)$. 

\end{proof}

\begin{assumption} \textbf{(Markovian structure)} In the sequel we will assume that the driver $f$ is Markovian, that is 
\[
	f(\omega, t, Z_t, U_t) = \bar{f} (X_t(\omega), Z_t, U_t)
\]
for some measurable $\bar{f}$. For convenience we simply write $f$ for $\bar{f}$. 
\label{markov}
\end{assumption}

\begin{corollary} Let $(Y^{\alpha,x,s}, Z^{\alpha,x,s}, U^{\alpha,x,s})$ be the unique bounded solution to the discounted BSDE
\begin{equation}
\begin{split}
	Y^{\alpha,x,s}_T &= Y^{\alpha,x,s}_t - \int_{t}^{T} (-\alpha Y^{\alpha,x,s}_{u} + f(X(t,s,x), Z^{\alpha,x,s}_u, U^{\alpha,x,s}_u)) du \\
	&+ \int_{t}^{T} (Z^{\alpha,x,s}_u)^*dW_u  
			   + \int_{t}^{T}\int_{B}U^{\alpha,x,s}_{s}(x)\tilde{N}(ds,dx),  
\end{split}
\label{Discounted_BSDE}
\end{equation}
on $[s,T]$ for some $s \geq 0$. We define a function $v^{\alpha}$ by $v^\alpha(s,x) = Y^{\alpha,x,s}_s$. Then, provided $v^{\alpha}$ is measurable, $Y^{\alpha,x,s}_t = v^\alpha(t,X(t,s,x))$ is a solution, and by uniqueness we also get that $v^\alpha(s,x)$ is bounded. It is also not hard to see that processes $Z$ and $U$ are Markovian, in other words the solution triplet $(Y_t,Z_t,U_t)$ can be represented as 
\[
	(v^\alpha(t,X_t), \xi^\alpha(t,X_t), \psi^\alpha(t,X_t))
\]
for some deterministic $v^\alpha, \xi^\alpha, \psi^\alpha$. For details on Markovian representations see, e.g. \cite{SCA}.
\label{mark_sol}
\end{corollary}

\noindent In what follows we will repeatedly use changes of measure to eliminate various parts of the driver in our BSDE. In view of Lemma \ref{Lemma_tmp}, in order to use the result of Theorem \ref{Estimate_coupl}, we need to ensure that under the new measure, the nonlinearity of the drift of the process $\{X_t\}_{t\geq 0}$ can be approximated as a uniform limit of Lipschitz functions. Hence we require the following assumption:

\begin{assumption} With the notation of Corollary \ref{mark_sol}, the function $\vartheta^{\alpha}: \bR^+ \times H \to \bR$, defined by
\[
	\vartheta^{\alpha}(t,x) := f(x,0,\psi^\alpha(t,x)) - f(x,0,0)
\] 
can be represented pointwise as a limit of a uniformly bounded family of Lipschitz (in $x$) functions for all $\alpha > 0$. 
\label{ALL}
\end{assumption}

\subsection{Ergodic BSDEs}

Now we use the same technique we employed in Theorem \ref{thm:discountedBSDEexist} to obtain a solution for the Ergodic BSDE
\begin{equation}\label{eq:EBSDE}
 Y_t = Y_T+\int_t^T[f(X^x_u, Z_u,U_{u})-\lambda] du - \int_t^T Z_u^* dW_u - \int_{t}^{T}\int_{B}U_{s}(x)\tilde{N}(ds,dx),  
\end{equation}
where $0\leq t\leq T< \infty,$ and $f:H\times H^{*}\times \bR \to \bR$ is a given function, $Y$ is a real-valued c\`adl\`ag stochastic process, $Z$ is a predictable process in $H^*$.
We change measure in such a way to get rid of the drift term, then take expectations, and then send $T$ to infinity. In our case, the generator depends on $\omega$ through the forward process $X(t,s,x)$, and we define measure $\bQ^{x,\alpha, T}$ to be such that the process
\[
\begin{split}
	\tilde{K}_t = \int_{t}^{T} (f(X(u,s,x), &Z^{\alpha_u,x,s}, U^{\alpha_u,x,s})-f(X(u,s,x), 0, 0)) du  \\
			& + \int_{t}^{T} (Z^\alpha_u)^*dW_u  + \int_{t}^{T}\int_{B}U^\alpha_{s}(x)\tilde{N}(ds,dx)
\end{split}
\]
is a $\bQ^{x,\alpha,T}$-martingale on $[s,T]$. Then, under $\bQ^{x,\alpha, T}$, we have 
\[
	e^{-\alpha s}v^\alpha(s,x) = E_{\bQ^{x,\alpha, T}}\Big[ e^{-\alpha T}v^\alpha(T,X(T,s,x)) + \int_{]s,T]} e^{-\alpha u} f(X(u,s,x), 0, 0)du\Big].
\]
As $|v^{\alpha}(t,X(t,s,x))|\leq C/\alpha$ for all $0 \leq t \leq T$, letting $T\to\infty$ we obtain
\[
	e^{-\alpha s}v^\alpha(s,x) = \lim_{T\to\infty}E_{\bQ^{x,\alpha, T}}\Big[\int_{]s,T]} e^{-\alpha u} f(X(u,s,x), 0, 0)du\Big].
\]
In order to proceed we notice that under the new measure our forward SDE takes the form
\begin{equation}
\begin{cases}
	dX_t = A(t)X_{t} + F^{\bQ}(t,X_{t})dt + G(t)dL_{t} \\
	X_{s} = x,	
\end{cases}
\label{Cauchy_final}
\end{equation}
where $F^{\bQ}(\cdot,\cdot)$ is the nonlinearity under measure $\bQ$ that includes new drift terms. 

\begin{lemma} The map $F^{\bQ}(t,x)$ is bounded and can be represented as a pointwise limit of a uniformly bounded family of Lipschitz functions.  
\end{lemma}
\noindent \begin{proof} We know explicitly the structure of $F^{\bQ}$. Define
\[
	\rho_s (x) := \begin{cases}
			   	\frac{f(x,\xi^{\alpha}(s,x),\psi^{\alpha}(s,x))-f(x,0,\psi^{\alpha}(s,x))}{||\xi^{\alpha}(s,x)||^2} (\xi^{\alpha}(s,x)) & \text{if } \xi^{\alpha}(s,x) \neq 0, \\
				0 & \text{otherwise.}	
			  \end{cases}
\]
where $Z_s = \xi^{\alpha}(s,X_s)$ and $U_s = \psi(s,X_s)$ as in Corollary \ref{mark_sol}. Using Lemma \ref{lm:gammaProcess}, define $\{\gamma_u\}_{u \geq s}$ to be such that 
\begin{equation}
	f(X(t,s,x),0,U_{t}) - f(X^x_t,0,0) =  \int_{B}U(v)\gamma_{t}(X(t,s,x),v)\nu(dv)
\label{gammaTmp}
\end{equation}
for all $t \geq 0$. Then $\bar{F}(t,X(t,s,x)) := F^{\bQ}(t,X(t,s,x)) - F_t(X(t,s,x))$ can be written as
\[
	\bar{F}(t,X(t,s,x)) =  G(t) \rho_t (X(t,s,x))+ \int_{B}  \gamma_t(X(t,s,x),r) \big[G(t) r\big]\nu(dr).
\]
The first argument is bounded due to the fact that $f$ is Lipschitz in $Z$. By arguments identical to Lemma 3.4 in \cite{Hu}, one can also show that it is a pointwise limit of Lipschitz functions. The second term depends on $X(t,s,x)$ only though the process $\gamma_{t}$. By (\ref{gammaTmp}) and Assumption \ref{ALL}, we conclude the result.

\end{proof}

\begin{lemma} For $v^\alpha$ defined as in Corollary \ref{mark_sol}, and for an arbitrary $x_0 \in H$, there exist bounds $C'$ and $C$ such that
\[
	|v^{\alpha}(s,x)-v^{\alpha}(s,x_0)|<C' (1 + \| x \|^2 + \| x_0 \|^2) \quad \text{and}  \quad \alpha |v^\alpha(s,x)| < C
\]
uniformly in $x$, $s$ and $\alpha$. 
\label{L1}
\end{lemma}
\noindent \begin{proof} With $\bQ^{x,\alpha, T}$ as above, we denote $\cP^\alpha(t,s) [f](x) = \bE^{\bQ^{x,\alpha, T}}[f(X(t,s,x))]$, where $X(t,s,x)$ is the mild solution to (\ref{Cauchy_final}). Then we obtain
\[
\begin{split}
	| v^\alpha(s,x) &- v^\alpha(t,x_0) | \\
						&\leq  e^{\alpha s}\bigg| \lim_{T\to\infty}E^{\bQ^{x,\alpha, T}}\Big[\int_{]s,T]} e^{-\alpha u} f(X(u,s,x), 0, 0)du\Big] \\
						& \qquad-  \lim_{T\to\infty}E^{\bQ^{x_0,\alpha, T}}\Big[\int_{]s,T]} e^{-\alpha (u-s)} f(X(u,s,x_0), 0, 0)du\Big] \bigg|  + C''|s-t|\\
						& \leq \int_s ^ {\infty} e^{-\alpha (u-s)} \big| \cP^\alpha(u,s) [f(\cdot,0,0)](x) - \cP^\alpha(u,s) [f(\cdot,0,0)](x_0) \big| du \\ & \qquad+ C''|s-t| \\
						 &\leq C' (1 + \| x \|^2 + \| x_0 \|^2) + C''|s-t|,
\end{split}
\]
for every $0 \leq t \leq s$, where $C'$ and $C''$ are independent of $\alpha$. For the last step we use the result of Theorem \ref{Estimate_coupl}. The inequality $\alpha |v^\alpha(s,x)| < C$ follows from Theorem \ref{thm:discountedBSDEexist}. 

\end{proof}

\begin{remark} We notice that, due to the periodic structure, $v^{\alpha}(t+T^*,x) = v^{\alpha}(t,x)$, and hence, for fixed $s \geq 0$ and $x_0 \in H$,
\[
\begin{split}
	\| v^{\alpha}(t,x) - v^{\alpha}(s,x_0) \| &\leq \| v^{\alpha} (t,x) - v^{\alpha}(t,x_0)\| + \| v^{\alpha}(t,x_0)  - v^{\alpha}(0,x_0)\| \\
				 &\leq C'(1 + \| x \|^2) + \sup_{u \in [s,s+T^*]} \| v^{\alpha}(u,x_0)  - v^{\alpha}(s,x_0)\| \\
				 & \leq C'(1 + \| x \|^2) + C''T^*.
\end{split}
\]
Given that $v^{\alpha}$ is uniformly Lipschitz in time, we have 
\[
	\| v^{\alpha}(t,x) - v^{\alpha}(s,x_0)\| \leq C'(1 + \| x \|^2)
\]
for some new constant $C'$.
\label{R1}
\end{remark}

\begin{lemma}There exists a bound $\bar{C}$, such that 
\[
	\| \nabla_x v^{\alpha}(t,x) \| \leq \bar{C}(1 + \| x \|^2).
\]
holds uniformly in $x,t$ and $\alpha$. 
\label{L2}
\end{lemma}
\noindent \begin{proof}
We begin by finding an estimate for the sensitivity of the process $\{X(t,\tau,x)\}_{t \geq \tau}$ with respect to the initial value $x$ (for the remainder of the proof we use the notation $X^{t,x}_s$ for $X(s,t,x)$). By standard arguments (see, e.g. \cite{Bismut_Levy}) one can show that, for any fixed $t \geq \tau$, there exists a constant $c_t > 0$ such that
\begin{equation}
	\bE \| \langle D_x X^{\tau,x}_t, h \rangle\|^2 \leq c_t \| h \|^2
\label{E2}
\end{equation}
holds for any direction $h \in H$. We also know that for any $t > 0$, there exists a probability measure $\bQ^{x,t} \sim \bP$, such that  
\[
	v^{\alpha}(t,x) = \bE^{Q,x,t} \bigg[ e^{-\alpha} v^{\alpha}(t+1,X^{t,x}_{t+1}) - \int_t^{t+1} e^{-\alpha (s-t)}f(X^{t,x}_s,0,0) ds \bigg].
\]
Since $ \nabla_x v^{\alpha}(t,x) = \nabla_x [ v^{\alpha}(t,x) - e^{-\alpha}v^{\alpha}(0,0)]$, we then obtain
\begin{equation}
\begin{split}
	\langle \nabla_x v^{\alpha}(t,x), h \rangle &=  e^{-\alpha}\bE^{Q,x,t}\bigg( \nabla_x \tilde{v}^{\alpha}(t+1,X^{t,x}_{t+1}) \langle D_x X(t+1,t,x),h \rangle \bigg) \\
						& \quad - \bE^{Q,x,t}\int_t^{t+1} e^{-\alpha(s-t)}\nabla _x f(X^{t,x}_s,0,0) \langle D_x X(s,t,x),h \rangle ds ,	
\end{split}
\end{equation}
where $\tilde{v}^{\alpha}(t,x) := v^{\alpha}(t,x) - v^{\alpha}(0,0)$. The last ingredient we need is the so called Bismut--Elworthy formula (for the L\'evy noise case see, e.g. \cite{Bismut_Levy}):
\[
	\bE \bigg[  \psi(X(s,t,x)) \int_t^s  G^{-1}(s)V^h_s dW_s \bigg] = (s-t) \langle h, D_x P(t,s) [\psi] (x) \rangle,
\]
where $V^h_s = \langle D_x X(s,t,x),h \rangle$, and $P(t,s)$ is the two parameter semigroup associated with $X(\cdot,t,x)$. Setting $\psi(\cdot) = \tilde{v}^{\alpha}(t+1,\cdot)$ and $\phi(\cdot) = f(\cdot,0,0)$, we notice that  
\[
	\langle D_x P(t,t+1) [\psi] (x) , h \rangle = \nabla_x \tilde{v}^{\alpha}(t+1,X^{t,x}_{t+1}) V^h_{t+1},
\]
and 
\[
	\langle D_x P(t,s) [\phi] (x) , h \rangle = \nabla_x f(X^{t,x}_s,0,0) V^h_s.	
\]
Therefore, using the Bismut--Elworthy formula twice, we have 
\[
\begin{split}
	\| \langle \nabla_x v^{\alpha}(t,x), h \rangle \|^2 &\leq 2 e^{-\alpha}  \bigg| \bE^{Q,x,t}   \big[ \tilde{v}^{\alpha}(t+1,X^{t,x}_{t+1}) \int_t^{t+1}  G^{-1}(s)V^h_s dW_s \big]  \bigg|^2\\
				& \quad + 2\bigg|\bE^{Q,x,t}\int_t^{t+1} e^{-\alpha(s-t)} \nabla _x f(X^{t,x}_s,0,0) V^h_s  ds \bigg|^2\\
				&\leq 2 \bE^{Q,x,t} \| \tilde{v}^{\alpha}(t+1,X^{t,x}_{t+1}) \|^2 \bE^{Q,x,t} \bigg( \int_t^{t+1}  \| G^{-1}(s)V^h_s \|^2 ds \bigg) \\
				&\quad + 2 \int_t^{t+1} \bigg( \frac{\bE^{Q,x,t} \| \phi(X^{t,x}_s) \|^2}{(s-t)^{-2}}  \bE^{Q,x,t} \int_t^s \| G^{-1}(u)V^h_u \|^2 du \bigg)ds.
\end{split}
\]
From Remark \ref{R1}, we know that 
\[
	\| \tilde{v}^{\alpha}(t+1,X^{t,x}_{t+1}) \|^2 \leq C'(1+\| X^{t,x}_{t+1} \|^2)^2, \quad \| \phi(\cdot) \| \leq C.
\]
The claim then follows taking into account (\ref{E2}) and the fact that $G^{-1}(t)$ is uniformly bounded.

\end{proof}

\begin{theorem} There exists a sequence $\alpha_n \to 0$, a bounded deterministic function $v:\bR^+ \times H \to \bR$ and a constant $\lambda\in\bR$, such that
\[
(v^{\alpha_n}(s,x)-v^{\alpha_n}(s,x_0)) \to v(s,x) \quad \text{ and } \quad \alpha_n v^{\alpha_n}(s,x) \to \lambda 
\]
for all $s \geq 0$, $x \in H$.
\label{EBSDE_construction}
\end{theorem}

\noindent \begin{proof} Since $H$ is a separable space, there exists a dense subset $V \subset \bR_{+} \times H $. On $V$ we can use a diagonal procedure to construct a sequence $\alpha_n \searrow 0$ such that 
\[
	(v^{\alpha_n}(s,x)-v^{\alpha_n}(s_0,x_0)) \to v(s,x) \quad \text{ and } \quad \alpha_n v^{\alpha_n}(s_0,x_0) \to \lambda
\]
for some function $v: V \to \bR$ and a real number $\lambda$. 
\noindent By Lemmas \ref{L1} and \ref{L2} we know that the functions $v^{\alpha}$ are locally Lipschitz in both time and space. We can therefore extend $v$ by continuity to the whole $\bR^+ \times H$, proving that 
\[
	v^{\alpha_n}(s,x) - v^{\alpha_n}(s,x_0) \to v(s,x)
\]
for all $x \in H$ and $s \geq 0$. We notice that, for $t \geq s$,
\[
\begin{split}
	\alpha_n v^{\alpha_n}(t,x) &= \alpha_n (v^{\alpha_n}(s_0,x_0)) + \alpha_n (v^{\alpha_n}(t,x) - v^{\alpha_n}(s_0,x_0)) \\
						&= \alpha_n v^{\alpha_n}(s_0,x_0) + \alpha_n (v^{\alpha_n}(t,x) - v^{\alpha_n}(t,X(t,s,x_0)) \\
						& \qquad +\alpha_n \int_{]s,t]} e^{-\alpha_n u} f(X(u,s,x_0),0,0)du \\
						&\to \lambda
\end{split}
\]
since
\[
\begin{split}
						\bigg| \alpha_n (v^{\alpha_n}(t,x) &- v^{\alpha_n}(t,X(t,s,x_0)) +\alpha_n \int_{]s,t]} e^{-\alpha_n u} f(X(u,s,x_0),0,0)du \bigg| \\
						&\leq  \alpha_n C'' |t-  s_0| + \alpha_n C' (1 + \| X(t,s_0,x_0) \|^2 + \| x_0 \|^2) \\
						&\to 0.
\end{split}
\]
We have thereby proven that $\lambda$ is indeed a constant independent of time. 

\end{proof}

\begin{theorem} Let $v$ and $\lambda$ be constructed as above. We also set $x_0 = 0 \in H$ and $s_0 = 0 \in \bR$ for the sake of simplicity. Then, if we define
\[
	Y^x_t = v(t,X^x_t),
\]
there exist processes $Z^x$ and $U^x$ such that the quadruple $(Y^x,Z^x,U^x,\lambda)$ solves the EBSDE
\[
 Y^x_t = Y^x_T+\int_t^T[f(X^x_{u}, Z^x_u,U^x_{u})-\lambda] du - \int_t^T (Z^x_u)^* dW_u - \int_{t}^{T}\int_{B}U^x_{s}(x)\tilde{N}(ds,dx)	
\label{EBSDE_final}
\]
for $0 \leq t \leq T < \infty$. Moreover, if there exists any other solution $(Y',Z',U',\lambda')$ that satisfies 
\begin{equation}
	|Y'_t| < c_x(1 + \| X^x_t \|^2),
\label{growth}
\end{equation}
for some constant $c$ that may depend on $x$, then $\lambda = \lambda'$.
\label{lambda_unique}
\end{theorem}
\noindent \begin{proof} We look at the discounted BSDE
\[
\begin{split}
	Y^{\alpha,x}_T &= Y^{\alpha,x}_t - \int_{t}^{T} (-\alpha Y^{\alpha,x}_{u} - \alpha v^\alpha(0,0) + f(X^x_t, Z^{\alpha,x}_u, U^{\alpha,x}_u)) du \\
	&\qquad + \int_{t}^{T} (Z^{\alpha,x}_u)^*dW_u  
			   + \int_{t}^{T}\int_{B}U^{\alpha,x}_{s}(x)\tilde{N}(ds,dx).  	
\end{split}
\] 
It is clear that the unique bounded solution is $Y^{\alpha,x}_t = v^\alpha(t,X^x_t)-v^\alpha(0,0)$. We remember that $|v^\alpha(s,X^x_s)-v^\alpha(0,0)| \leq (1 + \| X^x_t \|^2  )$. We conclude, by the dominated convergence theorem, that
\[
	\bE \int_0^T | Y^{\alpha,x}_t - Y^{\alpha_m,x}_t |^2dt \to 0 \quad \text{and} \quad \bE |Y^{\alpha,x}_T - Y^{\alpha_m,x}_T|^2 \to 0
\]
as $n \to \infty$. 

\noindent We now prove that the sequences $Z^{\alpha,x}$ and $U^{\alpha,x}$ are also Cauchy. Denote $\bar{Y} = Y^{\alpha_n,x} - Y^{\alpha_m,x}$, $\bar{Z} = Z^{\alpha_n,x} - Z^{\alpha_m,x}$, $\bar{U} = U^{\alpha_n,x} - U^{\alpha_m,x}$. We then have
\[
	\bar{Y}_T = \bar{Y}_t - \int_{t}^{T} \big(-\alpha \bar{Y}_{u} + \bar{f}(u)\big)du  + \int_{t}^{T} (\bar{Z}_u)^*dW_u  
			   + \int_{t}^{T}\int_{B}\bar{U}_{s}(x)\tilde{N}(ds,dx),
\] 
where $\bar{f}(u) = f(X^x_u, Z^{\alpha_n,x}_u, U^{\alpha_n,x}_u) - f(X^x_u, Z^{\alpha_m,x}_u, U^{\alpha_m,x}_u)$. By standard arguments, we know that, for any $\beta \geq 4K + 1/2$ (where $K$ is the Lipschitz constant of $f$), and $\beta > \max(\alpha_n,\alpha_m)$, we have
\[
\begin{split}
	e^{\beta t}\bE \| \bar{Y_t} \|^2 &+ \frac{1}{2}\int_t^T e^{\beta s} \bE \bigg( \| \bar{Z}_t \|^2 + \int_B \| \bar{U}_s(v) \|^2 \nu(dv) \bigg) dt \\
						     & \leq 	\bE \bigg[   \| \bar{Y}_T \|^2  + \frac{4}{2\beta -1} \int_t^T e^{\beta s} \| \delta f_s\|^2ds\bigg],
\end{split}
\]
where 
\[	
	\delta f_s = (\alpha_n - \alpha_m) Y^{\alpha_m,x}_s.
\]
By Theorem \ref{EBSDE_construction} and using the bound on $\bE [ \| X^x_t \|^2 ]$ obtained in Step 1 of Theorem \ref{Estimate_coupl}, we know that there exists $C = C(x)$ such that $\bE [\| Y^{\alpha_m,x}_s\|^2] \leq C$, and thus 
\[
	\bE \bigg[ \frac{4}{2\beta -1} \int_t^T e^{\beta s} \| \delta f_s\|^2ds \bigg] \leq \frac{4CT}{2\beta -1}(\alpha_n - \alpha_m)^2,
\]
and hence we immediately see that sequences $\{Z^{\alpha_n,x}\}_{n \geq 1}$ and $\{U^{\alpha_n,x}\}_{n \geq 1}$ are Cauchy. Denoting $Z^x$ and $U^x$ their corresponding limits, we get the first part of the result. 

\noindent In order to prove uniqueness, suppose there exists another solution $(Y', Z', U',\lambda')$ with polynomial growth. Let $\tilde Y = Y^x-Y'$, $\tilde Z=Z^x-Z'$, $\tilde{U} = U^x - U'$ and $\tilde \lambda=\lambda-\lambda'$. Then
\[
\begin{split}
	   \tilde Y_t = \tilde Y_T + \int_{]t,T]}  [f(X^x_u, Z^x_u,U^x_u) &- f(X^x_u, Z'_u,U'_u)- \tilde \lambda]du \\
	   &- \int_{]t,T]}  \tilde Z_u^*dW_u - \int_{]t,T]} \int_{B}\tilde{U}_{s}(r)\tilde{N}(ds,dr)
\end{split}
\]
By the standard Girsanov's argument there exists a probability measure $\bQ^T \sim \bP$ such that under $\bQ^T$ the process
\[
\begin{split}
	K_t = \int_{t}^{T} (f(X^x_u, Z^x_u, U^x_u)&-f(X^x_u, Z'_u, U'_u)) du \\
	&+ \int_{t}^{T} \delta \tilde{Z}_u^* dW_u  + \int_{t}^{T}\int_{B}\tilde{U}_{s}(r)\tilde{N}(ds,dr)
\end{split}
\]
is a martingale on $[0,T]$. Then we see that
\[
	\tilde \lambda = T^{-1} \bE^{\bQ^T}\big[\tilde Y_T - \tilde Y_0].
\]
Given the growth condition (\ref{growth}) and the estimate (\ref{X_estimate}), by sending $T \to \infty$ we obtain $\tilde{\lambda} = 0$ and thus the uniqueness of $\lambda$ is proven.

 \end{proof}

\noindent We are now ready to prove the main uniqueness result for Markovian solutions to our EBSDE, where by ``Markovian'' we mean that, if $Y$ is a solution, then there exists a continuous deterministic function $v$, such that $Y_t = v(t,X^x_t)$ for all $t > 0$. In the proof we will use the fact that the coefficients in the forward process are time dependant but $T^*$-periodic for some $T^* > 0$. Recalling the construction of the solution in Theorem \ref{EBSDE_construction} we immediately see that it is $T^*$-periodic in the first argument. Therefore, it is sensible to establish uniqueness in the class of Markovian solutions for which 
\begin{equation}
	v(t,x) = v(t+T^*,x) \quad \forall t > 0, x \in H.
\label{v_per}
\end{equation}

\begin{theorem} Let $(Y,Z,U,\lambda)$ and $(Y',Z',U',\lambda')$ be two Markovian solutions to the EBSDE (\ref{EBSDE_final}). If $Y,Y'$ satisfy the growth condition (\ref{growth}), $v,v'$ satisfy (\ref{v_per}) and $v'(0,0) = v(0,0)$, then $v = v'$ a.e.
\end{theorem}
\noindent \begin{proof}  From Theorem (\ref{lambda_unique}) we know that $\lambda = \lambda'$. We now show that in this case $Y = Y'$. Denoting $\tilde Y = Y^x-Y'$, $\tilde Z=Z^x-Z'$, $\tilde{U} = U^x - U'$ and defining $\bQ^T$ as in the proof of Theorem \ref{EBSDE_final}, we immediately have for all $t < T$
\[
	\tilde{Y}_t = \bE^{\bQ^T} [\tilde{Y}_T | \F_t]
\] 
for all $T$. Given the Markovian representation of our solutions we can rewrite the above as 
\begin{equation}
	\tilde{v}(t,x) = \bE^{\bQ^T} [ \tilde{v}(T,X^{t,x}_T) | \F_t ],
\label{v_tmp}
\end{equation}
where $\tilde{v}(t,x) := v(t,x) - v'(t,x)$. Now, since (\ref{v_tmp}) holds for any $T$, we obtain 
\[
	\tilde{v}(t,x) = \bE^{\bQ^{kT^*}} [ \tilde{v}(kT^*,X^{t,x}_{kT^*}) | \F_t ] = \bE^{\bQ^{kT^*}} [ \tilde{v}(0,X^{t,x}_{kT^*}) | \F_t ],
\]
for all $k$ such that $kT^* \geq t$. The next ingredient we will require is following estimate which can be shown with a technique identical to the one used to obtain (\ref{X_estimate}):
\[
	\bE^{\bQ^T} \big[ \| X^x_t \|^4 \big] < c (1 + \| x \|^4), \quad t \in [0,T]
\]
where $c$ is independent of $T$. We now notice that, for any $\epsilon > 0$, there exists $\delta > 0$ such that $|\tilde{v}(0,x)| \leq \epsilon$ if $\|x\| < \delta$, due to the fact that $\tilde{v}$ is locally Lipschitz $\tilde{v}(0,0) = 0$. Set $\tau = \inf \{ kT^*: \| X^{t,x}_{kT^*} \| < \epsilon, k \in \bN\}$. We then see that 
\[
\begin{split}
	|\tilde{v}(t,x)| &= \big| \bE^{\bQ^{kT^*}} [ \tilde{v}(0,X^{t,x}_{kT^* \wedge \tau}) | \F_t ] \big| \\
	&\leq \bE^{\bQ^{kT^*}}\big[ |\tilde{v}(0,X^{t,x}_{\tau})| \bone_{\{ \tau < kT^* \}}\big] + \bE^{\bQ^{kT^*}}\big[ |\tilde{v}(0,X^{t,x}_{kT^*})| \bone_{\{\tau \geq kT^* \}}\big] \\
				&\leq \epsilon + (\bQ^{kT^*}(\tau > kT^*))^{\frac{1}{2}} \big( \bE^{\bQ^{kT^*}} |\tilde{Y}_{kT^*}|^2  \big)^{\frac{1}{2}} \\
				&\leq \epsilon + C(\bQ^{kT^*}(\tau > kT^*))^{\frac{1}{2}} \big( \bE^{\bQ^{kT^*}} \big[ 1 + |X^{t,x}_{kT^*}|^4 \big] \big)^{\frac{1}{2}} \\
				& \to \epsilon.
\end{split}
\]
The last step of the derivation above follows from the fact that 
\[
	\bQ^{kT^*}(\tau > kT^*) \to 0 \text{ as } k \to \infty.
\]
In order to see this, we look at the discretised process $\{ X^{t,x}_{kT^*} \} _ {k \in \bN}$. We immediately see that it is irreducible. We therefore can prove the desired recurrence by following the proof of Theorem \ref{recurrence} with time step chosen as the first multiple of $T^*$ larger than $\tT$.

\end{proof}

\subsection{Alternative representation for $\lambda$}

In this section we show the representation of $\lambda$ as an integral with respect to a certain invariant measure. We established in Section 3.3 that there exists a unique invariant measure $\mu$ corresponding to the semilinear problem
\[
	dX_t = A(t)dt + F_t(X_t)dt + G(t)dL_t, \quad X_\tau = x.
\]
In particular, the following holds:
\[
	\int_{[0,T^*]\times H} P_s u(t,x) \nu(dt,dx) = \int_{[0,T^*] \times H}u(t,x)\nu(dt,dx),
\]
where $P_s$ is the corresponding semigroup. We recall that the Markovian solution to the EBSDE (\ref{EBSDE_final}) constructed in Theorem \ref{lambda_unique} is $T^*$-periodic, that is the quadruple $(Y,Z,U,\lambda)$ has a representation $(v,\xi,\psi,\lambda)$, where $v,\xi$ and $\psi$ are $T^*$-periodic in time. \begin{theorem} The value $\lambda$ in the EBSDE solution $(v,\xi,\psi,\lambda)$ satisfies
\[
	\lambda = \int_{[0,T^*] \times H} f(x,\xi(t,x),\psi(t,x)) \mu(dt,dx),
\]
where $\mu$ is the unique invariant measure.
\end{theorem}

\noindent \begin{proof} The invariance of $\mu$ implies that for any fixed times $T$ and $s \leq T$, and any bounded measurable function $u$ such that $u(t+T^*,x) = u(t,x)$ we have
\[
	\int_{[0,T^*]\times H} \bE u(T,X^{s,x}_T) \mu(dt,dx) = \int_{[0,T^*] \times H}u(t,x)\mu(dt,dx).	
\]
We write
\[
	v(t,x) = \bE^{\bP_{x,t}} \bigg[ v(T,X^{s,x}_T) + \int_t^T (f(X^{s,x}_t,\xi(s,X^{s,x}_t),\psi(s,X^{s,x}_t)) - \lambda) ds \bigg],
\]
where the subscript $(x,t)$ indicates that the forward equation was started at time $t$ with the value $x$. Then by the invariance property, integrating both sides with respect to $\mu$, we obtain the result. 
\end{proof}

\begin{remark} The representation above gives us an intuitive idea of how to interpret $\lambda$. If one thinks about the driver $f$ as a cost function of the optimally controlled dynamical system for the law of $X$, then $\lambda$ is the cost of one cycle.  
\end{remark}

\section{Applications}

\subsection{Classical Ergodic Control}
In this section we show how general ergodic control problems can be seen in the framework of EBSDEs for the case of controlled drift. Denote by $L: H \times \U \to \bR$ a bounded measurable cost function such that 
\[
	| L(x,u) - L(x',u)| \leq C \| x - x' \|,
\]
for some $C > 0$. We consider the problem of minimising 
\[
	J(x_0, u) = {\lim\sup}_{T\to\infty} T^{-1} \bE^{u,T}\bigg[\int_0^T L(X_t, u_t) dt\bigg],
\]
over the space $\U$ of controls, a separable metric space in which $u_t(\omega)$ takes values. We further assume that under $\bP^{u,T} \sim \bP$ the dynamics of the controlled process $X$ on $[0,T]$ are given by
\[
	dX_t = (A(t)X_t + F_t(X_t))dt + R(u_t)dt + \bigg( \int_B \gamma(u(t),y)\nu(dy)  \bigg) dt + G(t)dL_t, 
\] 
with $X_0 = x_0$. We further assume that $\| R(u) \| \leq C'$ and $\gamma(u(t),y)$ is a measurable function such that there exist a constant $0 \leq C < 1$ such that for every $u \in \U$ 
\[	
	-C(1 \wedge ||\xi||) \leq \gamma(u,\xi) \leq C(1 \wedge ||\xi||)
\]
for all $\xi \in B$. We define the Hamiltonian 
\begin{equation}
	f(x,z,r) = \inf_{u \in \U} \bigg\{ L(x,u) + zR(u) +\int_B \gamma(u,\xi)r(\xi) \nu(d\xi)\bigg\},
\label{infimum}
\end{equation}
where $x \in H$, $z \in H$ and $r : B \to \bR$. Immediately we notice that $f(x,0,0)$ is bounded. It is also easy to check that $f$ satisfies Assumptions \ref{ass:driver} and \ref{ass:gamma}. Therefore, the EBSDE with driver $f(x,z,r)$ admits a unique (in the class of processes with polynomial growth) Markovian solution $(Y,Z,U,\lambda)$. If the infimum in (\ref{infimum}) is attained, then, by a well known result (see \cite{SCA}), there exists (assuming the continuum hypothesis) a Borel-measurable function $\kappa: H \times H^* \times \L^{2}(B,\B,\nu,\bR) \to \U$ such that 
\[
	f(x,z,r) = L(x,\kappa(x,z,r)) + zR(\kappa(x,z,r)) + \int_B \gamma(\kappa(x,z,r),\xi)r(\xi)\nu(d\xi).
\]
\begin{theorem} Let the quadruple $(Y,Z,U,\lambda)$ be the unique Markovian solution satisfying $|Y_t| \leq c (1 + \| X_t \|^2)$ for all $t \geq 0$ and some $c > 0$. Then the following hold:
\begin{enumerate}[(i)]
\item For an arbitrary control $u \in \U$ we have $J(x_0,u) = \lambda$ if 
\[
	f(X_t,Z_t,U_t) = L(X_t,u(t)) + Z_t R(u(t)) + \int_B \gamma(u(t),\xi)r(\xi)\nu(d\xi) \quad d\bP \times dt - a.e.
\]
\item If the infimum is attained in (\ref{infimum}), then the control $\bar{u}(t) = \kappa(X_t,Z_t,U_t)$ verifies $J(x_0,\bar{u}) = \lambda$. 
\item Even if the infimum in (\ref{infimum}) is not attained, there exists a control $\{ \tilde{u}_t \}_{t \geq 0}$, such that $J(x_0,\tilde{u}) = \lambda$.
\end{enumerate}
\end{theorem}

\noindent \begin{proof} Identical to the proofs of Theorem 8 in \cite{Cohen_ergodic} and Theorem 5.1 in \cite{Andrew}.
\end{proof}

\subsection{Power plant evaluation}
In this section we present a model for power plant evaluation using Ergodic BSDEs. We show how due to the properties of gas and electricity the problem falls very naturally into the theoretical framework we have developed. We begin by defining a mathematical model of a power plant.
\begin{definition} We denote by $\{ E(t) \}_{t \geq 0}$ and $\{ G(t) \}_{t \geq 0}$ the electricity and gas price processes respectively. We assume that a power plant allows its owner to convert gas into electricity instantaneously, generating profit if $E(t) - cG(t) > 0$, where $c$ is some conversion constant. The quantity $X(t) := E(t) - cG(t)$ is called the spark spread. 
\end{definition}

\noindent In existing literature (for an overview see, for example, \cite{CCS}) the value of a power plant is approximated as a sum of spread options on spot power with different maturities, namely european options with payoffs $X_{T_j}^+$, where $\{ T_j, j \in J \}$ represent the future hours of production over the plant's lifetime. In other words 
\[
	VP_t = \sum_{j \in J}\exp(-r(T_j-t)) \bE^{\bQ} \bigg( (X_{T_j})^+ \bigg| \F_t \bigg).
\]
A flaw of this approach is that it relies heavily on the current state of the world, characterised by the short term dynamics of the electricity and gas prices. However, it is clear that one might want to evaluate the power plant before investing into its construction, and by the time the plant begins operation all the short term parameters will have changed. In the rest of the section we provide an alternative method for evaluation, assuming only that the price processes follow ergodic behaviour. In terms of the problem in question, this means that the present state is not important for the calculation of the long term (ergodic) average. 

We develop a slightly simplified model, where we do not give the dynamics of electricity and gas prices separately, but instead assume that the evolution of the spark spread $X$ is governed by the following equation: 
\begin{equation}
	d X_t = \theta_t (\kappa_t - X_t)dt + G(t) \big[ dW_t - \int_{B}x \tilde{N}(dt,dx) \big], \quad X_\tau = x,
\end{equation}
where $B = \bR \backslash \{0\}$, $\{ \theta_t \}_{t \geq 0}$ is a positive process that describes the rate of mean reversion, $\{ \kappa_t \} _ {t \geq 0}$ is a non-negative process of the mean and $\tilde{N}$ is a compensated Poisson random measure on $\bR_+ \times B$ with the compensator $\eta(dt,dx) =  \nu(dx) dt$. We also assume that all the processes are periodic in time with period $T^*=$ one year . The goal is to find the average yearly profit of the plant, namely 
\[
	\lambda = \lim_{T \to \infty} \frac{1}{T} \bE \int_t^T (X_s)^+ ds, 
\]
where $(x)^+ := \max(x,0)$. It is important to notice that, in reality, the difficulty in finding $\lambda$ comes from the fact that the vector of parameters $(\theta, \kappa, \nu)$ is not known exactly. Therefore, we face the risk averse problem of determining the worst-case average under a range of plausible parameters, namely 
\[
	\lambda = \inf_{u \in \U}\lim_{T \to \infty} \frac{1}{T}\bE^u \int_t^T (X^u_s)^+ ds, 
\]
where $\U$ denotes a space of possible values for $u = (\theta,\kappa,\nu)$, and under $\bP^u \sim \bP$ the dynamics of $X$ are given by
\[
\begin{split}
	d X_t = \theta_t (\kappa_t - X_t)dt &+ R(X_t,u(t))dt + \int_B \gamma(u(t),y)\nu(dy)dt \\
				&+ G(t) \big[ dW_t - \int_{B}x \tilde{N}(dt,dx) \big].
\end{split}
\]
These parameters control the rate of mean reversion through $R$ and the rate of spikes through $\gamma$. In order to make the model more realistic, without loss of clarity one can also consider the problem of minimising a generalised functional
\[
	\lambda = \inf_{u \in \U}\lim_{T \to \infty} \frac{1}{T}\bE^u \int_t^T L(X_s,u(s)) ds, 	
\]
where $L(x,u)$ incorporates a penalty corresponding to the perceived likelihood of the parameters being realised. Following exactly the same logic as in the derivation of  (\ref{infimum}), we define the Hamiltonian
\[
	f(x,z,r) = \inf_{u \in \U} \bigg\{ L(x,u) + zR(x,u) +\int_B \gamma(u,\xi)r(\xi)\nu(d\xi) \bigg\},	
\]
and proceed to solve the EBSDE with the driver $f$. 

\begin{remark} It is clear that once $\lambda$ is known, the risk-averse discounted expected revenue of the power plant with estimated lifetime of $N$ years can be calculated by  
\[
	v(N) = \lambda \int_0^N e ^ {-r(t)}dt,
\]
where $r(t)$ is a (deterministic) discount rate. 
\end{remark}

\begin{remark} As we mentioned at the beginning of this section, imposing the Ornstein--Uhlenbeck dynamics on the spark spread is restrictive. Ideally one would like to model electricity and gas processes separately. If we assume that the marginal price processes follow sums of OU processes (as in \cite{MBMM}, where the authors focus mainly on the copula-based approach) we end up with a two-dimensional problem, where the ergodicity required for the existence of a solution to EBSDE is obtained through the fact that the sum of ergodic processes is itself ergodic. The reason we present a simplified version is that it naturally demonstrates the theoretical framework we developed in previous chapters, and gives a clear illustration of how EBSDEs can be applied to this class of problems.
\end{remark}


\begin{thebibliography}{10}

\bibitem{Alru}
S.~Albeverio and B.~R\"udiger.
\newblock Stochastic integrals and the {L}\'evy--{It\^o} decomposition theorem
  on separable {B}anach spaces.
\newblock {\em Stochastic Analysis and applications}, 23:217--253, 2005.

\bibitem{Andrew}
A.~L. Allan and S.~N. Cohen.
\newblock Ergodic backward stochastic difference equations.
\newblock {\em arXiv:1509.00231}, 2015.

\bibitem{Barles}
G.~Barles, R.~Buckdahn, and E.~Pardoux.
\newblock Backward stochastic differential equations and integral-partial
  differential equations.
\newblock {\em Stochastics}, 60:57--83, 1997.

\bibitem{BL_impulse}
A.~Bensoussan and J.~L. Lions.
\newblock {\em Impulse control and quasi-variational inequalities}.
\newblock Gathier-Villars, Paris, 1982.

\bibitem{Bruss}
T.~Bruss.
\newblock A counterpart of the {B}orel--{C}antelli lemma.
\newblock {\em Journal of Applied Probability}, 17:1094--1101, 1980.

\bibitem{CCS}
Ren\'e Carmona, Michael Coulon, and Daniel Schwarz.
\newblock Electricity price modeling and asset valuation: a multi-fuel
  structural approach.
\newblock {\em Math Finan Econ}, 2013.

\bibitem{Cohen_gexp}
S.~N. Cohen.
\newblock Representing filtration consistent nonlinear expectations as
  \textit{g}-expectation in general probability spaces.
\newblock {\em Stochastic Processes and their applications}, 122(4):1601--1626,
  2012.

\bibitem{SCA}
S.~N. Cohen and R.~J. Elliott.
\newblock {\em Stochastic Calculus and Applications}.
\newblock Birkh\"auser, second edition, 2015.

\bibitem{Cohen_ergodic}
S.~N. Cohen and Y.~Hu.
\newblock Ergodic {BSDE}s driven by {M}arkov {C}hains.
\newblock {\em SIAM J. Control Optim.}, 51(5):4138--4168, 2013.

\bibitem{Coquet}
F.~Coquet, Y.~Hu, J.~Memin, and S.~Peng.
\newblock Filtration consistent nonlinear expectations and related
  \textit{g}-expectations.
\newblock {\em Probability theory and related fields}, 123(1):1--27, 2002.

\bibitem{Hu}
A.~Debussche, Y.~Hu, and G.~Tessitore.
\newblock Ergodic {BSDE}s under weak dissipative aassumptions.
\newblock {\em Stochastic Processes and their applications}, 121:407--426,
  2011.

\bibitem{Hu_Banach}
M.~Fuhrman, Y.~Hu, and G.~Tessitore.
\newblock Ergodic {BSDE}s and {O}ptimal {E}rgodic {C}ontrol in {B}anach
  {S}paces.
\newblock {\em SIAM Journal of Control and Optimization}, 48(3):1542--1566,
  2009.

\bibitem{Karoui_BSDE}
N.~El Karoui and L.~Mazliak.
\newblock {\em Backward Stochastic Differential Equations}.
\newblock Pitnam Research Notes, vol. 364, 1997.

\bibitem{Knable}
F.~Kn\"able.
\newblock Ornstein--{U}hlenbeck equations with time-dependent coefficients and
  {L}\'evy noise in finite and infinite dimensions.
\newblock {\em Journal of Evolution Equations}, 11:959--993, 2011.

\bibitem{Lindv}
T.~Lindvall.
\newblock {\em Lectures on the coupling methods}.
\newblock Dover Publications, 1992.

\bibitem{MBMM}
Th. Meyer-Brandis and M.~Morgan.
\newblock A dynamic {L}\'evy copula model for the spark spread.

\bibitem{Pardoux_Peng}
E.~Pardoux and S.Peng.
\newblock Adapted solution of a backward stochastic differential equation.
\newblock {\em Systems and Control Letters}, 14:55--61, 1990.

\bibitem{Zabczyk_Levy}
S.~Peszat and J.~Zabczyk.
\newblock {\em Stochastic partial differential equations with {L}\'evy noise}.
\newblock Cambridge University Press, 2007.

\bibitem{PRATO}
G.~Da Prato and J.~Zabzcyk.
\newblock {\em Stochastic Differential Equations in Infinite Dimensions}.
\newblock Cambridge University Press, 2008.

\bibitem{G.DaPrato1996}
G.Da Prato and J.~Zabczyk.
\newblock {\em Ergodicity for Infinite-Dimensional Systems}.
\newblock London Mathematical Society Lecture Notes, 1996.

\bibitem{PSXZ}
E.~Priola, A.~Shirikyan, L.~Xu, and J.~Zabczyk.
\newblock Exponential ergodicity and regularity for equations with {L}\'evy
  noise.
\newblock {\em http://arxiv.org/abs/1102.5553}, 2011.

\bibitem{PrZab}
E.~Priola and J.~Zabczyk.
\newblock Structural properties of semilinear {SPDE}s driven by cylindrical
  stable processes.
\newblock {\em Probability Theory and Related Fields}, 149:97--137, 2011.

\bibitem{Royer_jumps}
M.~Royer.
\newblock Backward stochastic differential equations with jumps and related
  non-linear expectations.
\newblock {\em Stochastic Processes and their applications}, 116:1358--1376,
  2006.

\bibitem{Bismut_Levy}
Bin Xie.
\newblock Uniqueness of invariant measures of infinite dimensional stochastic
  differental euqations driven by l\'evy noises.
\newblock {\em Potential Analysis}, 36(1):35--66, 2012.

\bibitem{XYZ}
J.~Yong and X.~Zhou.
\newblock {\em Stochastic Controls. Hamiltonian Systems and HJB Equations}.
\newblock Springer, 1999.

\end{thebibliography}

\end{document}